\def\Irr{{\rm Irr}}
\def\End{{\rm End}}
\def\Mat{{{\rm Mat}_n({\bf C})}}
\def\D{{\cal D}}
\magnification =\magstep1
\baselineskip =13pt
\overfullrule =0pt

\centerline {\bf HYPERGEOMETRIC FUNCTIONS ON REDUCTIVE GROUPS}
\vskip .7cm

\centerline {\bf M. Kapranov}

\vskip 1cm

An essential feature of the classical theory of hypergeometric functions
is the interplay between various ways of representing such functions:
by integrals, solutions of differential equations, or series. The theory
of $A$-hypergeometric functions [12-17] 
  unifies and generalizes a considerable 
part of the classical theory by putting into the forefront the concept
of a torus action. From a purely algebraic point of view, considering the
torus actions means that we are interested in how a polynomial breaks down
 into
a sum of monomials. Thus, an $A$-hypergeometric function is a function in
the coefficients of an indeterminate polynomial $f(x) = \sum_{\omega\in A} 
a_\omega x^\omega$
and $A\i {\bf Z}^m$ is a finite set of characters of the algebraic torus 
$({\bf C}^*)^m$,
i.e., of Laurent monomials.  For example, a typical Euler integral [16] is
 the function
$$J(f) = \int f(x)^{\alpha_0} x_1^{\alpha_1} ... x_m^{\alpha_m} dx_1...dx_m.
\leqno (0.1)$$
In the same vein, the $A$-hypergeometric series (or $\Gamma$-series)
[13-15] [17] are explicit power series in
the coefficients $a_\omega$. On the geometric side, 
the theory becomes allied with the algebro-geometric
theory of toric varieties [10] [18], the integral (0.1) being a period
 of a hypersurface
in such a variety. 

\vskip .2cm

The purpose of this paper is to give a generalization of the theory of
$A$-hypergeometric functions 
to the case when the torus is replaced by an arbitrary reductive
group $H$. The set $A$ is then a finite set of irreducible representations
 of  $H$.  The space of polynomials on a given set of monomials is replaced
 by the
space $M_A$ of functions on $H$ obtained as linear combinations of matrix
 elements of
representations from $A$. Thus a hypergeometric function is a function on 
the space $M_A$;
for instance, the Euler integral is written, similarly to (0.1), by 
$J(f) = \int f(x)^\alpha x^\sigma dx$ where $x^\sigma$ is a multivalued
 character of
$H$. The related algebro-geometric theory is that of spherical varieties, 
i.e., of
equivariant compactifications of reductive groups and their homogeneous
 spaces
coming from Gelfand pairs [3-5] [7-8] [32]. In fact, one of the reasons for
the author's thinking about this subject was to begin understanding  mirror
 symmetry
for hypersurfaces in spherical varieties and its relation to the Langlands
 duality for reductive groups.

The space $M_A$,
 being naturally a product of matrix spaces, can be regarded as a partial
compactification of another reductive group $G = \prod GL_{d(\omega)}$ where
$d(\omega), \omega\in A$ is the dimension of the representation. 
So the analogs of hypergeometric
series in this new situation are, again, series in matrix elements of 
irreducible
representations
of $G$. Their construction can be obtained  by generalizing the idea implicit
in [14], namely, expanding the delta-function along a subgroup $H\i G\subset 
M_A^*$
into a series
in matrix elements of representations of $G$ and then taking the 
termwise Fourier transform.
 This leads to a very general class of generalized power series.
Among their coefficients one finds the
 Clebsch-Gordan coefficients and other classical
quantities of representation theory. Given that in the toric case
the combinatorics of the coefficients of hypergeometric series is
at the basis of numerical predictions for the number of curves
in mirror symmetry [1] [25-27] [46], we give in \S 6 an example
where the coefficients involve the $3j$-symbols for the group $SU_2$.  

\vskip.2cm

The paper is organized as follows. In section 1 we introduce notation for
decomposing functions on a  reductive group into matrix elements
of irreducible representations, so as to emphasize the analogy with
 ordinary monomials
and polynomials and make writing formulas easier. Section 2 studies 
the spherical varieties $X_A, Y_A$
associated with a given set $A$ of irreducible representations.
These are not the most general spherical varieties, being
compactifications of reductive groups, of the kind  studied by De Concini 
and Procesi [7-8]. The varieties $Y_A$ are in fact, reductive algebraic
semigroups studied in [41] [49]. 
 In Section 3 
we discuss ``power series" on a group, i.e., series in matrix elements of
irreducible representations and their use in representing analytic functions
and distributions. Section 4 studies the Fourier transform of matrix
elements of irreducible representations of $GL(n)$ and the corresponding analog
of the gamma function. These questions have been investigated in [21-23]. 
Next, in Section 5 we introduce the $A$-hypergeometric system and
develop the analogs of all three approaches known for the toric case:
differential systems, Euler integrals, hypergeometric series.  
Finally, Section 6 is devoted to some examples.  

\vskip .2cm

I was recently informed by I.M. Gelfand that a part of the constructions
and results of this paper has been found by M.I. Graev (unpublished).

\vskip .2cm

My understanding of the subject of hypergeometric functions
owes very much to I.M. Gelfand and A.V. Zelevinsky, as it is I hope
clear from the content of the paper. 
The work on this paper started during the Taniguchi Symposium 
``Integrable Systems and Algebraic Geometry" held in Kobe and Kyoto
in July 1997. I am grateful to the Taniguchi Foundation and to the
organizers of the symposium for making this meeting possible.
I am also very much indebted to the referees of this paper
for pointing out numeruous shortcomings and inaccuracies in the
first version. Thanks are due to V. Ginzburg, B. Hall, B. Sturmfels
and V. Retakh
for bringing several relevant papers to my attention. 
This research was partly supported by an NSF grant
and by an A.P. Sloan Fellowship.

\vfill\eject

\centerline {\bf \S 1. Monomials and polynomials on a  group.}

\vskip 1cm

In this section we recall some well known properties of reductive groups.

\vskip .2cm

\noindent {\bf (1.1) Monomials.} Let $G$ be a reductive algebraic group
over {\bf C}. By ${\rm Irr}(G)$ we denote the set of isomorphism
classes of (finite-dimensional, algebraic) irreducible representations of $G$.
For $\alpha \in\Irr (G)$ we choose a representation $V_\alpha$
in the isomorphism class of $\alpha$, and denote by $\rho_\alpha: G\to
{\rm Aut}(V_\alpha)$ the $G$-action on $V_\alpha$. The value of $\rho_\alpha$
on $x\in G$ will be simply denoted by $x^\alpha$. We will denote by
$\alpha^-$ the label of the dual representation: $V_{\alpha^-} = V_\alpha^*$. 
We also set $d(\alpha) = \dim(V_\alpha)$.

\vskip .2cm

Let ${\bf C}[G]$ be the ring of regular functions on $G$, and $M_\alpha
\i {\bf C}[G]$ be the subspace spanned by the matrix elements of $x^\alpha$,
i.e., the space of functions of the form
$$f(x) = (C, x^\alpha) = {\rm tr}(C\cdot x^\alpha), \quad C\in {\rm 
End}(V_\alpha)^*
\simeq {\rm End}(V_\alpha). \leqno (1.1.1)$$
Here we identify ${\rm End}(V_\alpha)$ with its dual by means of the form
${\rm tr}(ab)$. The space $M_\alpha$ will be called the space
of monomials of type $\alpha$ on $G$. The group $G\times G$ acts on the left
on ${\bf C}[G]$ by the formula
$((g_1, g_2)f)(g) = f(g_1^{-1} g g_2)$, and it is well known that we
have a decomposition into irreducible subspaces
$${\bf C}[G] = \bigoplus_{\alpha\in\Irr (G)} M_\alpha ,\quad M_\alpha 
\simeq V_\alpha \otimes V_\alpha^*.$$
In other words, this means that any $f\in {\bf C}[G]$ can be uniquely written
as a finite sum of monomials:
$$f(x) = \sum_\alpha (C_\alpha, x^\alpha), \quad C_\alpha \in {\rm End}
(V_\alpha). \leqno (1.1.3)$$

\proclaim (1.1.4) Proposition. A function $f\in {\bf C}[G]$ is conjugacy
 invariant if and
only if all the coefficients in its expansion (1.1.3) are scalar matrices:
$C_\alpha  = c_\alpha \cdot {\rm Id}$. 

\vskip .1cm

The monomial corresponding to ${\rm Id}_{V_\alpha}$ is just the character of
 the representation,
which will be denoted by
$$s_\alpha(x) = ({\rm Id}, x^\alpha) = {\rm tr}(x^\alpha) . \leqno (1.1.5)$$

\vskip .2cm

\noindent {\bf (1.1.6) Examples.} (a) Let  $G=({\bf C}^*)^n$
be the algebraic torus of dimension $n$. Then $\Irr(G)$ consists of
1-dimensional representations forming a lattice ${\bf Z}^n$. For
$\alpha = (\alpha_1, ..., \alpha_n)\in {\bf Z}^n$ the corresponding monomial
is just the usual Laurent monomial
$$x^\alpha = x_1^{\alpha_1} ... x_n^{\alpha_n}, \quad x=(x_1, ..., x_n)
\in ({\bf C}^*)^n.$$

(b) Let $G=GL_n({\bf C})$. In this case, as well known,
irreducible representations of $G$ are labelled by sequences
$\alpha = (\alpha_1 \geq ... \geq \alpha_n)$ with $\alpha_i\in {\bf Z}$.
We will  sometimes call such sequences Young diagrams.
We will write $\alpha\geq 0$ if $\alpha_n\geq 0$. 
The case $\alpha = (1,0,...,0)$ corresponds to the standard representation
$V={\bf C}^n$ of $GL_n({\bf C})$. For arbitrary $\alpha$ the representation
space $V_\alpha$ will be denoted $\Sigma^\alpha (V)$, so that 
$\Sigma^\alpha$ is the Schur functor. The dimension of this space
is given by the formula
$$\dim(\Sigma^\alpha({\bf C}^n)) = d_n(\alpha):={\prod_{1\leq i<j\leq n}
 (\alpha_i-\alpha_j +j-i)\over \prod_{1\leq i<j\leq n} (j-i)}.\leqno (1.1.7)$$
The character $s_\alpha(x)$ is the Schur symmetric function
of the eigenvalues of $x$. Note that
$$\Sigma^{(\alpha_1+1, ..., \alpha_n+1)}(V) = \Sigma^{(\alpha_1, ..., 
\alpha_n)}
(V)\otimes \det(V), \quad x^{(\alpha_1+1, ..., \alpha_n+1)} = 
x^{\alpha}\cdot\det(x).
\leqno (1.1.8)$$
We set $|\alpha|=\sum \alpha_i$.

The entries of the  matrix  monomial $x^\alpha$ are homogeneous functions
in the matrix elements of $x$
of degree $|\alpha|$. They are 
polynomials
 if and only if $\alpha\geq 0$. The dual representation corresponds to
$\alpha^- = (-\alpha_n, ..., -\alpha_1)$.

\vskip .3cm

\noindent {\bf (1.2) Multiplication of monomials.} The problem of finding
 the product of two
monomials
$$(C_\alpha, x^\alpha) \cdot (C_\beta, x^\beta), \quad C_\alpha\in 
{\rm End}(V_\alpha),
C_\beta\in {\rm End}(V_\beta),$$
is reduced to the problem of decomposing $V_\alpha\otimes V_\beta$ into 
irreducibles.
Namely, let us write:
$$V_\alpha\otimes V_\beta \simeq \bigoplus_\gamma N_{\alpha\beta}^\gamma 
\otimes
V_\gamma, 
\quad N_{\alpha\beta}^\gamma = {\rm Hom}_G(V_\gamma, V_\alpha\otimes V_\beta),
\leqno (1.2.1)$$
so that $N_{\alpha\beta}^\gamma$ is the multiplicity space. Then we have
 the tensor product
of the coefficients
$$C_\alpha\otimes C_\beta \in {\rm End}(V_\alpha)\otimes {\rm End}(V_\beta) =
{\rm End}(V_\alpha\otimes V_\beta).$$ 
Let 
$$p= p_{\alpha\beta}: {\rm End}(V_\alpha\otimes V_\beta) = {\rm End}\biggl( 
\bigoplus_\gamma
N_{\alpha\beta}^\gamma \otimes V_\gamma\biggr) 
\to \bigoplus_\gamma {\rm End}(N_{\alpha\beta}^\gamma \otimes V_\gamma)$$
be the projection onto the diagonal blocks, and
$$ p_{\alpha\beta}^\gamma: {\rm End}(V_\alpha\otimes V_\beta) \to
 {\rm End}(N_{\alpha\beta}^\gamma)
\otimes {\rm End}(V_\gamma) \leqno (1.2.3)$$
be the $\gamma$th component of $p_{\alpha\beta}$. The following
is a well known general fact. 

\proclaim (1.2.4) Proposition. We have
$$(C_\alpha, x^\alpha)\cdot (C_\beta, x^\beta) = \sum_\gamma \bigl(
 {\rm tr}_1(p_{\alpha\beta}^\gamma(
C_\alpha\otimes C_\beta)), x^\gamma\bigr),$$
where 
$${\rm tr}_1: {\rm End}(N_{\alpha\beta}^\gamma) \otimes {\rm End}(V_\gamma) \to
{\rm End}(V_\gamma)$$
is the trace with respect to the second argument.

\vskip .3cm

\noindent{\bf (1.3) Multiplication and differentiation of monomials on
$GL_n$.} Consider the case $G=GL_n$. Then $x=\|x_{ij}\|$ is itself a monomial
corresponding to $\alpha = (1):= (1,0,...0)$, so $V_{(1)}=V$.
The decomposition of $V\otimes V_\alpha$ is given by the Pieri formula
(see, e.g.,  [11], Proposition 15.25):
$$V\otimes V_\alpha = \bigoplus_{\beta=\alpha+e_i} V_\beta,
\leqno (1.3.1)$$
where the sum is over $\beta = (\beta_1\geq ... \geq \beta_n)$
such that $\beta-\alpha$ is equal to the standard basis vector $e_i$
for some $i$. In particular, the decomposition is multiplicity-free,
so each occurrence of the trace (1.2.4) is taken over a 1-dimensional
space. Thus we can eliminate it from the explicit formulas,
saying that we have {\it canonical} maps
$p_{(1),\alpha}^\beta: {\rm End}(V)\otimes {\rm End}(V_\alpha)
\to {\rm End}(V_\beta)$, even though there is no canonical
element in the 1-dimensional space of maps $V\otimes V_\alpha\to V_\beta$. 

We can write therefore:
$$(a,x)\cdot (b_\alpha, x^\alpha) = \sum_{\beta=\alpha + e_i} 
\biggl( p_{(1), \alpha}^\beta (a\otimes b_\alpha), \, x^\beta\biggr).
\leqno (1.3.2)$$
Let now $\partial = \|{\partial\over\partial x_{ij}}\|$. This is a 
matrix-valued
vector field on $G$; more precisely, $\partial \in {\rm Vect}(G)\otimes
{\rm End}(V^*)$. If $a$ is an $n$ by $n$ matrix, we denote by
$(a,\partial)$ the scalar differential operator 
$\sum a_{ji} \partial/\partial x_{ij}$. The following fact is
 the generalization of the formula
$d(x^\alpha)/dx = \alpha x^{\alpha-1}$  for scalar monomials.

\proclaim (1.3.3) Proposition. For $a\in {\rm End}(V^*)$ and $b_\alpha \in
{\rm End}(V_\alpha)$ we have
$$(a,\partial)\cdot (b_\alpha, x^\alpha) = 
\sum_{\beta=\alpha -e_i} (\alpha_i+n-i) \cdot \biggl( (p_{(1), \beta}^\alpha)^t
(a\otimes b_\alpha),\, x^\beta\biggr),$$
where 
$$(p_{(1), \beta}^\alpha)^t: {\rm End}(V^*\otimes V_\alpha)\to {\rm 
End}(V_\beta)$$
is the transpose of
$$p_{(1), \beta}^\alpha: {\rm End}(V)\otimes {\rm End}(V_\beta)
\to {\rm End}(V_\alpha).$$

\noindent {\sl Proof:} Let $L$ be the vector space of translation
invariant vector fields on $\Mat$. It is an irreducible
representation of $G\times G$, isomorphic to $V^*\otimes V$ while
$M_\alpha \simeq V_\alpha\otimes V_\alpha^*$. The LHS and the RHS
 of the proposed
equality are $G\times G$-equivariant maps
$$l, r: L\otimes M_\alpha \to \bigoplus_{\beta = \alpha - e_i} M_\beta.$$
The source of $l,r$ is the tensor product of two irreducible
representations. The target is a direct sum of several distinct irreducibles.
Further, in the decomposition of the source into irreducibles
each irreducible enters with multiplicity $\leq 1$ by Pieri's formula.
Thus we have the following fact.

\proclaim (1.3.4) Lemma. Suppose there is a vector $v\in L\otimes M_\alpha$
whose projection to each $M_\beta$-isotypic component from the target of
$l,r$,
is nonzero and such that $l(v)=r(v)$. Then $l=r$ on the entire
$L\otimes M_\alpha$. 

Take now $v={\rm Id}\otimes {\rm Id}\in L\otimes M_\alpha$. We claim
that $l(v)=r(v)$.   
Indeed, this statement reduces to a  computation involving only the characters
$s_\alpha$. Consider $s_\alpha$ as a symmetric function of $n$ variables
$t_1, ..., t_n$, and let $D=\sum \partial/\partial t_i$. Our statement 
that $l(v)=r(v)$ thus
 reduces to the identity
$$(Ds_\alpha)(t) = \sum_{\beta = \alpha - e_i} (\alpha_i+n-i) s_\beta (t). 
\leqno
(1.3.5)$$
This identity follows at once from the Weyl character formula [11]
once we notice the denominator of that formula
is annihilated by $D$ and thus
can be treated as a constant. 

Further, in the case when $\alpha_i+n-i\neq 0$ for all $i$,
the equality (1.3.5) shows also that the projection of $v$ to
any $M_\beta$-isotypic component we are interested in, is nonzero,
so Proposition 1.3.3 is true for such an $\alpha$. Let us now identify
$M_\alpha$ with $M_{\alpha_1+s, ..., \alpha_n+s}$ for any $s\in {\bf Z}$
in the standard way. Then, we get two  families of maps
$$l(s), r(s): L\otimes M_\alpha \to \bigoplus_{\beta = \alpha - e_i} M_\beta,
\quad s\in {\bf Z}$$
which obviously depend on $s$ in a polynomial way. By the above,
$l(s)=r(s)$ for all except possibly finitely many $s$. Therefore
$l(s)=r(s)$ for any $s$ and thus $l=r$ for any $\alpha$. 
Proposition is proved.

\vfill\eject

\centerline {\bf \S 2. Some algebraic geometry related to monomials and 
polynomials.}

\vskip 1cm

\noindent {\bf (2.1) The space $M_A$. The (in)homogeneity condition.}
We start with a  reductive group $H$ and a finite set $A\i \Irr(H)$ of 
irreducible representations.
Let
$$M_A = \bigoplus_{\omega\in A} M_\omega  = \biggl\{ \sum_{\omega\in A}
 (a_\omega, x^\omega), 
\,\,\, a{_\omega} \in {\rm End}(V_\omega)\biggr\} \i {\bf C}[H] \leqno
 (2.1.1)$$
be the space of polynomials on monomials from $A$. So $M_A$ is a 
finite-dimensional
subspace  in ${\bf C}[H]$ invariant with respect to both left and right
 $H$-actions. The coefficients
$a_\omega$ serve as matrix coordinates in $M_A$, identifying it with 
$\bigoplus_{\omega\in A}
{\rm End}(V_\omega)$. We can thus identify $M_A$ with $M_A^*$ via the
 form
$\sum_{\omega\in A} {\rm tr}(a_\omega b_\omega)$ and view both $M_A, M_A^*$
 as partial
compactifications of the group $G = \prod_{\omega\in A} GL(V_\omega)$.
 We will be interested in generic behavior of functions from $M_A$. Let
$$\rho_A = \bigoplus_{\omega\in A} \rho_\omega:
 H \to \prod_{\omega\in A} GL(V_\omega) = G \i M_A^*-\{0\}
\leqno (2.1.2)$$ be the direct sum of representations from $A$. We will
 assume that $\rho_A$ is injective
so that we can regard $H$ as a subgroup of $G$. 

In the sequel we will impose one of the following conditions on $\rho_A(H) \i 
G$
(or, equivalently, on $A\i \Irr(H)$). Clearly, one of these conditions is
always satisfied.

\proclaim (2.1.3) Homogeneity condition. The image $\rho_A(H)$ contains the
 subgroup of scalars
$${\bf C}^* = \bigl\{ (\lambda, ..., \lambda) \in \prod GL(V_\omega), \,\, 
\lambda\in
{\bf C}^*\bigr\}.$$

\proclaim (2.1.4) Inhomogeneity condition. The intersection $\rho_A(H)\cap 
{\bf C}^*$ is finite.

We shall say that the pair $(H,A)$ is homogeneous, resp. inhomogeneous,
if it satisfies (2.1.3), resp. (2.1.4).

Given $(H,A)$ satisfying (2.1.4), we can construct $(\bar H, \bar A)$ 
satisfying (2.1.3)
so that $M_{\bar A} = M_A$. namely, let $\bar H$ be the image of ${\bf 
C}^*\times H$
under the multiplication map $$p: {\bf C}^*\times H \to G = \prod 
GL(V_\omega),
\quad (\lambda, x)\mapsto \lambda x. \leqno (2.1.5)$$
Then each $V_\omega$, $\omega\in A$, is an irreducible $\bar H$-module. 
Denoting
$\bar\omega$ its isomorphism class, we get a bijection
$$A\to\bar A, \,\, \omega\mapsto \bar\omega, \,\,\, V_{\bar\omega} = 
V_\omega.$$

\vskip.1cm

\noindent {\bf (2.1.6) Examples. }(a) Let $H= GL_n$. A set $A\i \Irr(GL_n)$ 
satisfies the homogeneity 
condition if and only if for any $\omega = (\omega_1 \geq ... \geq 
\omega_n)\in A$
the sum $|\omega| = 
\sum \omega_i$ is the same. 

(b) let $H=SL_n$. Then any set $A\i \Irr(SL_n)$ satisfies the inhomogeneity 
condition.

\vskip .3cm

\noindent {\bf (2.2) The varieties $Y_A, X_A$.} We assume that $A$ satisfies 
the
homogeneity condition. Let $Y_A$ be the Zariski closure of $\rho_A(H)$ in 
$M_A^*$. This is
a conic variety, and we denote by $X_A = P(Y_A)\i P(M_A^*)$ its 
projectivization. The
following properties of $X_A, Y_A$ are obvious from the construction.

\proclaim (2.2.1) Proposition. (a) $Y_A$ is invariant under the $H\times 
H$-action in
$M_A^*$ and contains $\rho_A(H)$ as an open orbit. \hfill\break
(b) $X_A$ is projective, invariant under the action of $P(H)\times P(H)$ and 
contains
$P(\rho_A(H))$ as an open orbit. \hfill\break
(c) Functions from $M_A$ are precisely the restrictions to $\rho_A(H)$ of 
linear functions
on $Y_A\i M_A^*$. \hfill\break
(d) ${\bf C}[Y_A]$, the algebra of regular functions on $Y_A$, is the 
subalgebra of
${\bf C}[H]$ generated by $M_A$.

Call a subset $S\i \Irr(H)$ {\it monoidal} if $0$ 
(the label for the trivial representation) is in $S$ and for each $\alpha, 
\beta\in S$ and
each embedding $V_\gamma\i V_\alpha\otimes V_\beta$
we have $\gamma\in S$. Let $\langle A\rangle$ be the minimal
monoidal subset generated by $A$.

\proclaim (2.2.2) Proposition. 
The subalgebra ${\bf C}[Y_A]$ has the form
$${\bf C}[Y_A] = \bigoplus_{\alpha\in \langle A\rangle} M_\alpha \i 
\bigoplus_{\alpha\in \Irr(H)} M_\alpha = {\bf C}[H], \quad M_\alpha
=V_\alpha\otimes V_\alpha^*.$$

\noindent {\sl Proof:} ${\bf C}[H]$, as an $H\times H$-module, splits
into a direct sum of distinct irreducible representations $M_\alpha$
for all $\alpha\in\Irr(H)$. Thus ${\bf C}[Y_A]$, being an
 $H\times H$-submodule, should have the form $\bigoplus_{\alpha\in S} M_\alpha$
where $S\i \Irr(H)$ is some subset. Since ${\bf C}[Y_A]$ is also an algebra
generated by the $M_\omega, \omega\in A$, we find that $S=\langle A\rangle$.

\proclaim (2.2.3) Corollary. The variety $Y_A$ has a natural structure
of an algebraic semigroup containing the group $H$.

\proclaim (2.2.4) Proposition. The number of $H\times H$-orbits on $Y_A$
is finite.

\noindent {\sl Proof:} This is a consequence of the fact that
each irreducible representation of $H\times H$ enters ${\bf C}[Y_A]$
no more than once, see [44].

\proclaim (2.2.5) Definition. If $(H,A)$ does not satisfy the homogeneity 
condition, then we define $X_A:=X_{\bar A}$,
$Y_A = Y_{\bar A}$, where $(\bar H, \bar A)$ is constructed in (2.1.5). 

\vskip .3cm

\noindent {\bf (2.3) Tannakian point of view on $Y_A$.}
We can generalize the construction of $Y_A$ as follows. 
Let $S$ be any monoidal subset in $\Irr(H)$
generated by a finite subset.  Then the subspace
$M[S]:= \bigoplus_{\alpha\in S} M_\alpha$ is a finitely generated
subalgebra in ${\bf C}[H]$.
The corresponding affine algebraic variety will be denoted
by $Y_{[S]}={\rm Spec} \, M[S]$. Thus $Y_A$ is obtained when
$S=\langle A\rangle$. 

To the set $S$ we associate the category ${\cal R}={\cal R}_S$ 
 all regular representations $V$ of $H$
such that each irreducible component of $V$ is isomorphic to some $V_\alpha$,
$\alpha\in S$. This is a monoidal category,
i.e., if $V,W\in {\cal R}$, then $V\otimes W\in {\cal R}$.
Denote by $f: {\cal R}\to {\rm Vect}$ the forgetful functor to the
monoidal category of vector spaces (with the usual tensor product).
Thus $f(V)$ is $V$ regarded as a vector space. The following fact
is a version of the well known Tannaka-Krein theorem reconstructing
a group from the category of its representations.

\proclaim (2.3.1) Proposition. (a) $Y_{[S]}$ is an affine algebraic
semigroup equipped with a semigroup homomorphism $i_S: H\to Y_{[S]}$
with an open dense image. The map $i_S$ is injective if and only if
every representation $V_\alpha, \alpha\in\Irr(H)$ can be embedded
into $V_\beta^* \otimes V_\gamma$ for $\beta, \gamma\in S$.
\hfill\break
(b) The correspondence $S\mapsto Y_{[S]}$ establishes a bijection
between finitely generated  monoidal subsets $S\i \Irr(H)$ and 
isomorphism classes
of homomorphisms $i: H\to Y$ where $Y$ is an affine semigroup
and ${\rm Im}(H)$ is open dense. 
\hfill\break
(c) Points of $Y_{[S]}$ are in bijection with monoidal
natural transformations $a$ of functors from $f$ to itself, i.e., with
systems of endomorphisms $a_V: V\to V$ given for all $V\in {\cal R}$
and satisfying the following properties:\hfill\break
(i) If $\phi: V\to W$ is a morphism of representations, then
 $\phi a_V = a_W\phi$;\hfill\break
(ii) $a_{V\oplus W} = a_V\oplus a_W$, $a_{V\otimes W} = a_V\otimes 
a_W$.\hfill\break
The semigroup structure on $Y_{[S]}$ is given by the composition of the 
operators:
$(ab)_V = a_V \circ b_V$. 

 The varieties $Y_{[S]}$ are precisely the
reductive algebraic semigroups studied in [41] [49].

\vskip .1cm

\noindent {\bf (2.3.2) Example.} Suppose that $H$ is a torus. 
Then $\Irr(H)$ is the lattice of characters of $H$,
a monoidal subset in it is just a sub-semigroup, $M[S] = {\bf C}[S]$
is the semigroup ring and Proposition 2.3.1 is the classification
of affine toric varieties.

\vskip .2cm

\noindent {\bf (2.4) The Newton polytope and the orbit structure of $X_A$.}
The varieties $X_A, Y_A$ belong to the following general class of varieties 
with group
action [3] [33].

\proclaim (2.4.1) Definition. (a) Let $K$ be a reductive group. An algebraic 
subgroup
$\Delta\i K$ is called spherical if some Borel subgroup in $K$ has a dense 
orbit
in $K/\Delta$.  \hfill\break
(b) An algebraic variety $M$ with $K$-action is called spherical if it
 contains an open orbit isomorphic to $K/\Delta$, where $\Delta$ is a 
spherical subgroup.
\hfill\break
(c) For any affine algebraic group $G$ its rank ${\rm rk}(G)$ is defined
as the dimension of any maximal torus. The rank of a spherical
$K$-variety $M$ is defined by ${\rm rk}(M) = {\rm rk}(K) - {\rm rk}(\Delta)$
where $\Delta$ is as in (b). 

An equivalent formulation of the condition in (a) is that each irreducible
 representation
of $K$ occurs  in the space ${\bf C}[K/\Delta]$ not more than once. 
Pairs $(K, \Delta)$ with this property are known as Gelfand pairs.

The example relevant for us is $K = H\times H$, and
 $\Delta = \{(x,x)\}\simeq H$ being
the diagonal. Clearly, $Y_A$ is spherical with this choice of $K, \Delta$.
 Similarly, $X_A$ is spherical.

Usually, one includes normality in the definition of spherical varieties.
Since this is unnatural in our context, we will not do this
and will take care in applying results which are stated in the literature
for normal varieties.

\vskip .1cm

Choose a maximal torus $T\i H$. Let $\Lambda = {\rm Ch}(T)$ be its character 
lattice
(i.e., the lattice of weights of $H$), and set $\Lambda_{\bf R} = \Lambda 
\otimes {\bf R}$.
Let also $W$ be the Weyl group of $H$ (with respect to $T$). 
We will identify $\Irr(H)$ with the set of dominant weights in $\Lambda$.

Given a representation
$\omega\in\Irr(H)$, we denote by $Q_\omega \i \Lambda_{\bf R}$
the convex hull of all the weights of $T$ on $V_\omega$. 
As well known, $Q_\omega = {\rm Conv}(W\cdot \omega)$.
 For a finite subset $A\i \Irr(H)$
as before, we denote by $Q=Q_A$ the convex hull of the union of the
 $Q{_\omega}$ for $\omega\in A$. 
This is a $W$-invariant convex polytope in $\Lambda_{\bf R}$. It will be 
called the
{\it Newton polytope} of $A$. 

\proclaim (2.4.2) Theorem. (a) $H\times H$-orbits on $Y_A$ are in 
order-preserving bijection with $W$-orbits on
 faces of $Q_A$ (including the empty face). Denote by $Y(\Gamma)$ the orbit
 corresponding to a face
$\Gamma$. \hfill\break
(b) Therefore the orbits on $X_A$ are the projectivizations $X(\Gamma):=
P(Y(\Gamma))$ for all nonempty faces $\Gamma\i Q_A$ (considered
up to $W$-action).
\hfill\break
(c) Explicitly, $Y(\Gamma)$ is the orbit of the vector $e(\Gamma)
\in\prod_{\omega\in A} {\rm End}(V_\omega)$ whose components
$e(\Gamma)_\omega: V_\omega\to V_\omega$ are defined as follows.
Write the weight decomposition $V_\omega = \bigoplus_{\lambda\in\Lambda}
V_\omega^\lambda$ with respect to the $T$-action. Then $e(\Gamma)_\omega$ is 
the block-diagonal operator with respect to this decomposition
whose diagonal block on $V_\omega^\lambda$ is zero if $\lambda\notin\Gamma$
and is the identity if $\lambda\in\Gamma$. 
\hfill\break
(d) Each $Y(\Gamma)$ is a spherical variety of rank equal to
${\rm rk}(H)-{\rm codim}(\Gamma)$.\hfill\break
(e)  In particular, minimal (closed) orbits correspond to $W$-orbits on the 
vertices
 of $Q_A$. For such a vertex $\lambda$ the orbit $X(\lambda)$
is isomorphic to 
$(H/P_\lambda)\times (H/P_\lambda)$, where $P_\lambda \i H$ is the parabolic 
subgroup whose
relative Weyl group is the stabilizer $W_\lambda \i W$.

\noindent {\bf (2.4.3) Examples.} (a) Let $H=GL_n = GL(V)$ and 
let $A$ consist of one element $\omega = (1)$, so that $V_\omega = V$
is the tautological representation. Then $Y_A = \Mat$.
 Orbits
of $H\times H$ on $\Mat$ are just sets of matrices with fixed rank. 
Thus there are $n+1$ orbits. The polytope $Q_A$ is the $(n-1)$-dimensional
simplex, the convex hull of the $n$ basis vectors in $ \Lambda_{\bf R}
= {\bf R}^n$, and the Weyl group acts by permutations of the vertices.
Thus there is only one orbit on faces of each dimension, and there are
$n+1$ orbits.

\vskip .1cm

(b) Let $H$ be a semisimple group of adjoint type and let
$A =\{\omega\}$ consist of one element which is strictly dominant. 
The variety $X_{\bar A}$ corresponding to the homogenized pair
$(\bar H, \bar A)$ from (2.1), is
the minimal compactification of $H$ constructed by De Concini and Procesi [7].
In this case $Q$ is the ``permutohedron" ${\rm Conv} (W\cdot\omega)$, 
with $|W|$ vertices forming one
orbit. As shown in [7], $X_{\bar A}$ has in this case one closed orbit,
 isomorphic to $(H/B)\times (H/B)$,
where $B$ is a Borel subgroup in $H$. 

\vskip .3cm

\noindent {\bf (2.5) Proof of Theorem 2.4.2.} The statement is covered
by the existing general theory of spherical varieties [3] [8] [33] and
especially that of reductive algebraic semigroups [41] [49]. 
Because of various notational differences, however,
it seems
the easiest to give a more self-contained proof, using the ideas of 
the cited papers. 

\proclaim (2.5.1) Lemma. Let 
$$K=H\times H,\quad  L=T\times T,\quad  V=\bigoplus_{\omega\in A} M_\omega, 
\quad v=(1,...,1).$$
Then any $K$-orbit in $Y_A = \overline{K\cdot v}$ intersects the closure of
$L\cdot v$.

\noindent {\sl Proof:} (cf. [41], \S 3, Lemma 3.)
Let $H((t))$, resp. $H[[t]]$ be the group of ${\bf C}((t))$,
resp. ${\bf C}[[t]]$-points of $H$, and similarly for $T$. 
The Iwahori theorem  (also known as the Cartan decomposition)
 for reductive groups over local fields such as
${\bf C}((t))$ gives
$$H((t)) = H[[t]] \cdot T((t))\cdot H[[t]].\leqno (2.5.2)$$
If $x\in Y_A$, we have $x=\lim_{t\to 0} g(t) v$ for some meromorphic
analytic curve $g(t)$, $0< |t| < \epsilon$ in $H$. Such a curve can be viewed
as an element of $H((t))$. Factoring now $g(t) = f_1(t) \gamma(t) f_2(t)$
with $f_i(t)\in H[[t]]$ and $\gamma(t)\in T((t))$, we get, taking into
account that $v$ is the unit element of $H$,
$$x=\lim_{t\to 0} g(t) v =\lim_{t\to 0}
(f_1(t) \gamma(t) f_2(t)) v = f_1(0) \biggl(\lim_{t\to 0}\gamma(t) v\biggr)
f_2(0),$$
which means that  $\lim_{t\to 0} \gamma(t) v$ exists in $Y_A$ and that
its
 $H\times H$-orbit  contains
$x$, proving the lemma. 

\vskip .1cm

Let ${\rm Ch}(L)$ be the lattice of characters of $L$
and ${\rm Ch}_{\bf R}(L) = {\rm Ch}(L)\otimes {\bf R}$. 
Let now $V$ be any algebraic representation of the torus $L$ and $V=\bigoplus
_{\lambda\in {\rm Ch}(L)} V^\lambda$ be its weight decomposition.
For a vector $v\in V$ denote by $v_\lambda\in V^\lambda$ its
component with respect to this decomposition. Set
$$A_L(v) = \bigl\{ \lambda: v_\lambda\neq 0\bigr\}\i 
{\rm Ch}(L),\quad {\rm Wt}_L(v) = {\rm Conv} (A_L(v))
\i {\rm Ch}_{\bf R}(L).\leqno (2.5.3)$$
Thus ${\rm Wt}_L(v)$ is a convex polytope
known as the weight polytope of $v$, see [18] Ch.5, \S 1. 
For any nonzero vector $w\in V$ let $P(w)\in P(V)$ be
the corresponding point of the projective space. 
The following fact is well known, see, e.g., {\it loc. cit.} Proposition 1.9:

\proclaim (2.5.4) Proposition. Let $\alpha: {\bf C}\to L$ be any rational
map and let $v'\in V$ be any nonzero vector
representing the point $\lim_{\tau\to 0} P(\alpha(\tau)\cdot v)\in P(V)$.
Then there is a face $\Gamma\i {\rm Wt}_L(v)$ such that $v'_\lambda = 0$
for $\lambda\notin A_L(v)\cap \Gamma$ while
$v'_\lambda$ is a nonzero multiple of $v_\lambda\neq 0$ for $\lambda\in
A_L(v)\cap\Gamma$. 

 Setting now as before $V=\bigoplus_{\omega\in A} M_\omega$, we see that 
 $A_L(v)$ is the set of all the weights in all the representations
$V_\omega$, embedded diagonally into $\Lambda\oplus\Lambda = {\rm Ch}(L)$.
So ${\rm Wt}_L(v)=Q_A$. For any weight $\lambda\in A_L(v)$
the component $v_\lambda$ is the identity map of the
corresponding weight space. Thus Proposition 2.5.4 implies that
any point in the closure of $P(L\cdot v)\i X_A$  is in the $L$-orbit
of the point described in Theorem 2.4.2 (c). Together with Lemma
2.5.1, this proves parts (a)-(c) of Theorem 2.4.2. 

Let us show (e). Let   $\lambda$ be a vertex of $Q_A$;
by applying the Weyl group transformations, if necessary, we can assume that
$\lambda$ is dominant. Since $\lambda$ is a vertex,
the corresponding weight subspace of $\bigoplus_{\omega\in A} V_\omega$
is 1-dimensional, and thus the vector $e(\lambda)$ of operators
has only one nonzero component, namely $e(\lambda)_\lambda: V_\lambda\to 
V_\lambda$.
 Further, this component is just the identity operator
on the 1-dimension weight $\lambda$ subspace extended by 0 along the
other weight spaces. In other words,  
 $e(\lambda)_\lambda = \xi_\lambda\otimes \xi_\lambda^*\in V_\lambda\otimes
V_\lambda^*$, where $\xi_\lambda$ is the highest vector in $V_\lambda$
and $\xi_\lambda^*$ is the highest vector in $V_\lambda^*$.
This implies (e).

Finally, part (d) follows from Theorem 7.3 of [33] applied to
the normalization of $Y_A$ once we translate the contravariant
description of spherical varieties (via fans) used there
into the dual covariant description (via polytopes) used here.

\vskip .3cm

\noindent {\bf (2.6) Orbit structure of semigroups.}
For future reference we note a slight generalization of Theorem
2.4.2. Let $S\i\Irr(H)$ be any finitely generated monoidal
subset and $Y_{[S]} = {\rm Spec}(M[S])$ be the corresponding affine
semigroup from (2.3). The group $H\times H$ acts on $Y_{[S]}$ via the
homomorphism $i_S: H\to Y_{[S]}$ and has an open orbit whose
stabilizer contains $H$. So $Y_{[S]}$ is an affine spherical
variety. Let $T, \Lambda$ be as in (2.4) and let $C_S\i\Lambda_{\bf R}$ be
the convex hull of all the $T$-weights of all the representations
$V_\alpha, \alpha\in S$. This is a polyhedral cone. The Weyl group
$W$ naturally acts on $C_S$. The proof of the next theorem is similar
to that of 2.4.2.

\proclaim (2.6.1) Theorem. $H\times H$-orbits on $Y_{[S]}$ are in bijection
 with
$W$-orbits on nonempty faces of the cone $C_S$. The orbit $Y(\Gamma)$
corresponding to a face $\Gamma$ is a spherical variety of rank equal to
$\dim(\Gamma)$. 

We can also exhibit a distinguished point 
 $e(\Gamma)\in Y(\Gamma)$. For this, we use the description of points of
$Y_{[S]}$ from Proposition 2.3.1(c). Let $\cal R$ be
 the monoidal category of representations corresponding to $S$. For every
$V\in {\cal R}$ we have the decomposition $V= V'_\Gamma \oplus V''_\Gamma$
where $V'_\Gamma$ is the direct sum of all weight subspaces
whose weights lie in $\Gamma$ and $V''_\Gamma$ is the direct sum
of all the other weight subspaces. To describe a point
$e(\Gamma)$ we will desribe the corresponding system of
operators $e(\Gamma)_V: V\to V$, $V\in{\cal R}$.

\proclaim (2.6.2) Proposition. 
The collection of operators $e(\Gamma)_V = 
{\rm Id}_{V'_\Gamma}\oplus 0_{V''_\Gamma}$ satisfies the conditions
of Proposition 2.3.1 (c) and thus defines a point $e(\Gamma)\in Y_{[S]}$.
The orbit $Y(\Gamma)$ contains $e(\Gamma)$.

\vskip .2cm

\noindent {\bf (2.7) The degree of $X_A$.} In the notation of (2.4),
 let $\Lambda^\vee = {\rm OP}(T)$ be the lattice of 1-parameter subgroups in 
$T$,
i.e., the dual lattice to $\Lambda$. Denote by $R\i\Lambda$ the root system of 
$H$ and
choose
a system of positive roots $R_+\i R$. Let also $R^\vee\i\Lambda^\vee$ be the 
system of
 coroots [45]
and $R^\vee_+$ be the system of positive coroots. Recall that $R$ and $R^\vee$ 
are in
bijection; for $\alpha\in R$ let $\alpha^\vee\in R^\vee$ be the corresponding 
coroot.
Each $\alpha^\vee$ is a linear function on the real space $\Lambda_{\bf R}$. 

\proclaim (2.7.1) Theorem. The degree of the projective variety
 $X_A\i P(M_A^*)$ is equal to
$${1\over \prod(d_i-1)!} \int_{Q_A} \prod_{\alpha\in R_+} \langle \alpha^\vee, 
\lambda\rangle d\lambda,$$
where $d\lambda$ is the Lebesgue measure on $\Lambda_{\bf R}$ normalized
 so that ${\rm Vol}
(\Lambda_{\bf R}/\Lambda)=1$, and the $d_i$ are the characteristic exponents 
of $H$, i.e., the
degrees of the polynomial generators of the algebra of $W$-invariant
 polynomials on
$\Lambda_{\bf R}$. 

This theorem was proved by B. Kazarnovski [31] in 1987. A more general result
for arbitrary spherical varieties was proved by M. Brion [4] in 1989.
Note that the function $\prod_{\alpha\in R_+} \langle \alpha^\vee, 
\lambda\rangle$ is
the famous Weyl volume element for the Langlands dual group $G^L$, see [45]. 
The maximal torus $T^L$ of $G^L$ is dual to $T$, so the quotient of $T$ by its 
maximal compact
subgroup is canonically identified with $\Lambda_{\bf R}$. Thus the integral in
 (2.7.1) can be viewed as the volume
of a conjugacy invariant domain in $G^L$ whose intersection with $T^L$ is the 
preimage of
$Q_A$ under the map $T^L\to T^L/T^L_c = \Lambda_{\bf R}$.

\vskip .3cm

\noindent {\bf (2.8) The $A$-discriminant.} Let $A\i \Irr(H)$ be an arbitrary
finite subset (satisfying or not the homogeneity condition).
 As in the toric case [18],
we set
$$\nabla_A^0 = \biggl\{ f\in M_A| \,\,\, \exists \, x_0\in H: \, f(x_0) = 
d_{x_0}f = 0\biggr\}$$
and let $\nabla_A\i M_A$ be the Zariski closure of $\nabla_A^0$. We call 
$\nabla_A$ the
$A$-discriminantal variety; in the case when it has codimenion 1, the equation 
of
$\nabla_A$ is called the $A$-discriminant and denoted by  $\Delta_A$.
Thus $\Delta_A$ is a polynomial function on $M_A = \bigoplus_{\omega\in A} 
{\rm End}(V_\omega)$.

\proclaim (2.8.1) Proposition. (a) $\nabla_A$ is conic and is invariant under
 both left and right actions of $H$ in $M_A$.\hfill\break
(b) The projectivization of $\nabla_A$ is projectively dual to $X_A$.

Extending this construction, for any nonempty face $\Gamma\i Q_A$ we denote
by $\nabla_{A, \Gamma}\i M_A$ the conic variety whose
projectivization is projectively dual to the orbit $X(\Gamma)\i X_A$.
This variety is irreducible, since $X(\Gamma)$ is.
 It is clear that $\nabla_{A,\Gamma}$
is also invariant under the left and right actions of $H$ on $M_A$. 
In the case when $\nabla_{A, \Gamma}$ is a hypersurface, we will denote
by $\Delta_{A, \Gamma}$ the irreducible equation of $\nabla_{A, \Gamma}$
(defined up to a constant factor).

\vskip .2cm

In the case when $H$ is a torus, the discriminants have been studied
in [18] from the point of view of monomials  entering into their expansion.
In the present case it is natural to view the space
$M_A$ on which the $\Delta_{A, \Gamma}$ are defined, as a partial
compactification of the group $G= \prod GL(V_\omega)$ and expand
the discriminants not in ordinary monomials but in the matrix elements
 of irreducible representations
of $G$: 
$$\Delta_{A,\Gamma} \biggl( \sum_{\omega\in A} (a_\omega, x^\omega)\biggl) =
 \sum_{\alpha\in\Irr (G)} (c_\alpha,
{\bf a}^\alpha).\leqno (2.8.2)$$
Here ${\bf a} = (a(\omega))_{\omega\in A}$, thus $\alpha = 
(\alpha(\omega))_{\omega\in A}$
with each $\alpha(\omega)$ being an element of $\Irr (GL(V_\omega))$, i.e., a 
Young diagram.
The coefficient $c_\alpha$ lies in
${\rm End}\biggl( \bigotimes_{\omega\in A} \Sigma^{\alpha(\omega)} 
(V_\omega)\biggr)$.
The condition of left and right $H$-invariance of $\nabla_{A,\Gamma}$ implies 
the following
quasi-homogeneity property of its equation.

\proclaim (2.8.3) Corollary. Suppose that
$\Delta_{A, \Gamma}$ is defined. There is a 1-dimensional character $\chi = 
\chi_\Gamma$
 of $H$ such that each
coefficient $c_\alpha$ of $\Delta_{A,\Gamma}$ satisfies
$$c_\alpha \in {\rm End}\biggl(\biggl( \bigotimes_{\omega\in A} 
\Sigma^{\alpha(\omega)} V_\omega
\biggr)^\chi\biggr),$$
where the superscript $\chi$ means the subspace of $\chi$-quasiinvariants.

\vskip .2cm





\vskip 2cm 

\centerline {\bf \S 3. Power series and distributions on a group.}

\vskip 1cm

\noindent {\bf (3.0) Orthogonality relations.} 
 Let $G$ be a reductive
group, as in \S 1, and $G_c\i G$ a compact form. Denote by $d^*x$
the regular
invariant volume form on $G$ normalized so that $\int_{G_c} d^*x=1$.
The orthogonality relations for the matrix elements of irreducible
representations have the form:
$$\int_{G_c} (a_\alpha, x^\alpha) \cdot (b_\beta, x^\beta) d^*x =
\cases{ 0, & if $\alpha\neq \beta^-$; \cr
d(\alpha)^{-1}{\rm tr}(a_\alpha b_\alpha^t), & if $\alpha =\beta^-$.\cr}
\leqno (3.0.1)$$
Here $d(\alpha)=\dim(V_\alpha)$ and $V_{\alpha^-} = V_\alpha^*$.
In particular, for the characters we have:
$$\int_{G_c} s_\alpha(x) s_\beta(x) d^*x = \cases{0&if $\alpha\neq \beta^-$;\cr
1&if $\alpha =\beta^-$.\cr} \leqno (3.0.2)$$

\vskip .2cm

\noindent {\bf (3.1) Power series and analytic functions.} 
By a power series on $G$ we mean a series of the form
$$f(x) = \sum_{\alpha\in\Irr(G)}(a_\alpha, x^\alpha), \quad a_\alpha
\in {\rm End}(V_\alpha). \leqno (3.1.1)$$
We are interested in representing analytic functions on (some domains in) $G$
in this way.  Representing of functions on $G_c$ in this way is the content
of the Peter-Weyl theorem which, together with the orthogonality
relations (3.0.1) gives the following.

\proclaim (3.1.2) Proposition. Let $f(x)$ be an analytic function in some 
neighborhood
of $G_c$ in $G$. Then $f$ can be written as the sum of an absolutely convergent
series (3.1.1) with
$$a_\alpha = d(\alpha) \int_{G_c} f(x) x^{\alpha^-} d^*x.$$

\noindent {\bf (3.2) Example. The case of $GL_n({\bf C})$.} For $G=GL_n$
we take $G_c = U_n$, the unitary group. In this case, as well known [29]
$$d^*x = {1\over C} {dx\over \det(x)^n}, \quad dx = \prod dx_{ij}, 
\quad C= \int_{U_n} {dx\over \det(x)^n} = 
{(2\pi)^{n(n+1)/2}\over \prod_{i=1}^{n-1} i!}.
\leqno (3.2.1)$$
Also, if $\alpha = (\alpha_1, ..., \alpha_n)$ is a Young diagram, then
$\alpha^- = (-\alpha_n, ..., -\alpha_1)$. So the formula for the
coefficients becomes
$$a_\alpha = {d(\alpha)\over C} \int_{U_n} f(x) x^{(-n-\alpha_n, ..., 
-n-\alpha_1)} dx.
\leqno (3.2.2)$$
Recall that $U_n$ is the Shilov boundary of the domain
$$D_n = \bigl\{ x\in {\rm Mat}_n({\bf C}): I-xx^* >0\bigr\}. \leqno (3.2.3)$$
Here, as usual, the notation $> 0$ means that the (Hermitian) matrix in
question is positive definite. So we have the following fact.

\proclaim (3.2.4) Proposition. Let $f$ be a function analytic on a neighborhood
of $U_n$ in $GL_n({\bf C})$. If $f$ can be continued to an analytic
function in a connected domain  in ${\rm Mat}_n({\bf C})$  containing 0,
then $f$ can be analytically continued to the entire $D_n$ and, in addition,
the coefficients $a_\alpha$ in (3.2.2) are zero unless $\alpha\geq 0$.

The reader can consult [35] for the general study of domains of convergence of
power series on reductive groups. 




\vskip .2cm

\noindent {\bf (3.3) Rings of formal series and convergence.} 
Denote by ${\bf C}[[G]]$ the vector space of all formal power series of the
form (3.1.1). This space is not a ring. However, we have the following fact.

\proclaim (3.3.1) Proposition. Let ${\cal D}(G)$ be the ring of
all regular differential operators on $G$. Then ${\bf C}[[G]]$ is naturally
a ${\cal D}(G)$-module, containing the space ${\bf C}[G]$ of regular functions.

\noindent {\sl Proof:} The structure of ${\bf C}[G]$-module on
${\bf C}[[G]]$ follows from Proposition 1.2.4 describing the
product of two monomials. If $L$ is a left-invariant vector
field on $G$, then $L$ preserves each $M_\alpha$ and so
acts on ${\bf C}[[G]]$. The ring ${\cal D}(G)$ is generated by ${\bf C}[G]$
and invariant vector fields. The check that the relations
holding in ${\cal D}(G)$ are satisfied, is standard.

\vskip .2cm

Let now $S\i\Irr(G)$ be a finitely generated monoidal subset.
Let also $T\i G$ be a maximal torus and $\Lambda$ its character
lattice, and $C_S\i \Lambda_{\bf R}$ be the cone defined in (2.6).
We shall say that $S$ is {\it strictly convex} if $C_S$ is a strictly
convex cone, i.e., it does not contain linear subspaces and this
$0\in\Lambda_{\bf R}$ is a vertex of $C_S$. Let $Y_{[S]}$ be the affine
semigroup associated to $S$.
 Theorem 2.6.1 implies that $Y_{[S]}$ contains a unique $G\times G$-fixed point
which we also denote $0$. 

Let $M[[S]]$ be the completion of the ring $M[S] = {\bf C}[Y_{[S]}]$
at $0$. This completion can be viewed as the ring of formal series
$$f(x) = \sum_{\alpha\in S} (a_\alpha, x^\alpha), \quad a_\alpha\in {\rm End}
(V_\alpha).\leqno (3.3.1)$$
Because of the strict convexity, the multiplication of two such
series is now a well defined, purely algebraic operation accomplished
via Proposition 1.2.4. Here is a useful sufficient condition for
convergence of a series from $M[[S]]$.

\proclaim (3.3.2) Theorem. Let $f\in M[[S]]$ be a formal series.
Suppose that there is a holonomic system of linear PDE on $G$ with
regular singularities (i.e., a holonomic regular ${\cal D}(G)$-module
[2]) of which $f$ is a solution. Then $f$ converges to an analytic function
defined on a neighborhood of 0 in $Y_{[S]}$. In particular, the
domain of convergence of $f$ in $G$ is nonempty. 

\noindent {\sl Proof:} For a smooth algebraic variety $Z$ let ${\cal D}_Z$
be the sheaf of rings of differential operators on $Z$. Let $\cal L$ be
the sheaf of ${\cal D}_G$-modules corresponding to the holonomic system
in question. Let also $G'\i Y_{[S]}$ be the open $G\times G$-orbit, 
i.e., the set of all invertible elements of $Y_{[S]}$ as a semigroup.
Then we have a surjective morphism of reductive groups $\sigma: G\to G'$.
By Hironaka's theorem, there exists a resolution of singularities 
$p: Y\to Y_{[S]}$ bijective over $G'$ and such that $p^{-1}(0)$ is a divisor
with normal crossings. Then $p^*f$ is a function on the formal neighborhood
of $p^{-1}(0)$ in $Y$. In particular, for any $y\in p^{-1}(0)$ we have the
germ $(p^*f)_y$ lying in $\widehat{\cal O}_{Y,y} \simeq {\bf C}[[y_1, ..., 
y_n]]$, the completion of the local ring of $Y$ at $y$. Let $\phi$ be the
composition
$$G\buildrel\sigma\over\to G'\buildrel p^{-1}\over\to p^{-1}(G')
\hookrightarrow Y.$$
By the general theory, see [2], we have that $\phi_*{\cal L}$ is
a holonomic regular ${\cal D}_Y$-module. It is also clear that
each $(p^*f)_y$, $y\in p^{-1}(0)$, is a formal solution of this
module, i.e.,
$$(p^*f)_y\in {\rm Hom}_{{\cal D}_Y} (\phi_*{\cal L}, \widehat{\cal 
O}_{Y,y}).$$
Now, it is the main property of holonomic regular systems, see, e.g., 
[2], Prop. 14.8, that any formal solution is convergent. Thus each
$(p^*f)_y$ is convergent in some neighborhood of $y$. This means
that there is a neighborhood $U\supset p^{-1}(0)$ in $Y$ such that
$p^{*}(f)$ converges to an analytic function in $U$. But
$U\cap p^{-1}(G')$ is identified with an open set in $G'$, so we get
that $f$ converges to an analytic function on a nonempty open subset in
$G'$ (and hence  to a function on a domain of  $G$). Theorem is
proved.

\vskip .3cm

\noindent {\bf (3.4) Distributions on a group.} As with usual power series,
series of the form (3.1.1) can be used to represent distributions on the group
as well. We need the following examples which follow at once from the
Peter-Weyl theorem.

\vskip .2cm

\noindent {\bf (3.4.1) Delta-function at 1.} It follows from
 Proposition 3.1.2 that the series
$$\delta_1(x) = \sum_{\alpha\in\Irr(G)} (d(\alpha)\cdot {\rm Id}, x^\alpha) = 
\sum_{\alpha\in\Irr(G)}
d(\alpha) {\rm tr}(x^\alpha)$$
(sum over all the representations) represents the Dirac delta-function 
situated at the
unity of $G$. This is an equality of distributions on $G_c$.

Accordingly, the delta-function situated at any $x_0\in G$ can be
written as
$$\delta_1(x_0^{-1}x) = \delta_1(xx_0^{-1}) = \sum_{\alpha\in\Irr(G)} 
d(\alpha)\biggl(
(x_0^{-1})^\alpha, x^\alpha\biggr).$$

\vskip .2cm

\noindent {\bf (3.4.2) Delta-function along a subgroup.}  Let $H\i G$ be a 
reductive
subgroup. For any $\alpha\in\Irr(G)$ the space $V_\alpha$ splits canonically
into direct sum
$$V_\alpha = V_\alpha^H \oplus V_\alpha',$$
where $V_\alpha^H$ is the space of $H$-invariants, and $V_\alpha'$ is the sum
of the other isotypical components. Thus we have the embedding:
$${\rm End}(V_\alpha^H) \i {\rm End}(V_\alpha).$$
Let $I_\alpha(H)$ be the image of the identity operator under this embedding.
Then:
$$\delta_H(x) = \sum_{\alpha\in\Irr(G)} d(\alpha) \cdot(I_\alpha(H), 
x^\alpha)$$
is the delta-function along $H$. Again,
this is an equality of distributions on $G_c$. Note that 
$$\delta_H(x) = \int_{y\in H_c} \delta_1(xy^{-1})    d^*y,$$
where $H_c = H\cap G_c \i H$ is a compact form. This follows from the fact 
that the
operator $\int_{y\in H_c} \rho_\alpha(y) d^*y$ is equal to the projection onto
$V_\alpha^H$ along $V_\alpha'$.

\vskip 2cm 

\centerline {\bf \S 4. The matrix gamma function.}

\vskip 1cm

\noindent {\bf (4.1) Motivation.} In the one-dimensional Fourier
analysis it is well known that the Fourier transform of a monomial
$x^\alpha$, $\alpha\in {\bf C}$, is again a monomial
$x^{-1-\alpha}$ but taken with the coefficient $\Gamma(1+\alpha)$,
the Euler gamma function. For application to hypergeometric series
we need a generalization of this to the case when $x$
is an $n\times n$ matrix and $\alpha$ is a Young diagram,
i.e., we need to study the Fourier transform of the matrix
elements of irreducible representations of $GL_n({\bf C})$.
This generalization was essentially done in [21-23],
and we just formulate it in the form needed for applications. 

\vskip .3cm

\noindent {\bf (4.1) Multivalued monomials.} Throughout this section we set:
$$ V={\bf C}^n,\,\, G=GL_n({\bf C}) = {\rm Aut}(V),\,\,
E=\Mat = {\rm End}(V).$$
We denote ${\bf g} = gl_n({\bf C})$ the Lie algebra of $G$. 
Let also $j: G\hookrightarrow E$ be the embedding. We also set
$$\Irr = \Irr(G) = \bigl\{ \alpha = (\alpha_1 \geq ... \geq \alpha_n) \in {\bf 
Z}^n\bigr\}.$$
For $\alpha\in\Irr$ we denote $V_\alpha = \Sigma^\alpha V$ the space
of the corresponding representation. 

By means of the standard basis in $V={\bf C}^n$ we identify
 $\Lambda^n V\simeq {\bf C}$. Accordingly, for any $s\in {\bf Z}$
and $\beta\in \Irr$ we identify $V_{(\beta_1+s, ..., \beta_n+s)}$
with $V_\beta$ as a vector space, so that the $G$-action differs by
$$\rho_{(\beta_1+s, ..., \beta_n+s)}(x) = \rho_\beta(x) \cdot \det(x)^s.$$
Let now
$${\bf I} = \biggl\{ \alpha = (\alpha_1, ..., \alpha_n)\in {\bf C}^n\,
\bigl| \alpha_i-\alpha_{i+1}\in {\bf Z}_+\biggr\}. \leqno (4.1.1)$$
Every $\alpha\in {\bf I}$ has the form
$$\alpha = \beta + (s, ..., s), \quad s\in {\bf C}, \beta\in \Irr.
\leqno (4.1.2)$$
The set {\bf I} serves to label
 multivalued holomorphic representations of $G$. Explicitly, if $\alpha\in 
{\bf I}$ has the
 form (4.1.2), we set $V_\alpha = V_\beta$ and define the
$G$-action by $\rho_\alpha(x) = \rho_\beta(x) \cdot\det(x)^s$.
As usual, we write $x^\alpha$ for $\rho_\alpha(x)$. This is a multivalued
representation, but the corresponding representation of the Lie algebra
{\bf g} is well-defined and irreducible. 

Let ${\cal M}_\alpha$ be the local system on $G$ spanned by all
 the determinations (branches) of all the matrix elements of $x^\alpha$. 
More precisely, if $\alpha$ is as in (4.1.2), and ${\cal L}_s$ is the
1-dimensional local system spanned by all the branches of $\det(x)^s$,
then ${\cal M}_\alpha = M_\beta\otimes {\cal L}_s$, with $M_\beta=
{\rm End}(V_\beta)$. Note that by construction ${\cal M}_\alpha$
is a subsheaf in ${\cal O}_G$, the sheaf of holomorphic functions. 

\proclaim (4.1.3) Proposition. (a) ${\cal M}_\alpha$ is a
 $G\times G$-equivariant local system on $G$.
As such, it is irreducible. \hfill\break
(b) We have an isomorphism of  equivariant local systems ${\cal M}_\alpha^* 
\simeq
{\cal M}_{(-\alpha_n, ..., -\alpha_1)}$.

\proclaim (4.1.4) Corollary. There is a natural action of the Lie algebra 
${\bf g}\oplus {\bf g}$
on each stalk of ${\cal M}_\alpha$ and thus on any space $\Gamma(U, {\cal 
M}_\alpha)$
where $U$ is a connected and simply connected open set in $G$. This action
is irreducible, of type $V_\alpha^*\otimes V_\alpha$. 

\vskip .2cm

\noindent {\bf (4.2) Fourier integrals over Hermitian matrices.}
Denote by 
 ${\rm Herm}$ the subspace of Hermitian matrices in $E$. 
So Herm is an $n^2$-dimensional totally real subspace
in $E$. Let ${\rm Herm}_+$ be the cone of positive definite
matrices in Herm, and let ${\cal H}\i E$ be the Siegel upperhalf plane
$${\cal H} = {\rm Herm} + i {\rm Herm}_+ = \{x: i(x^*-x) \geq 0\}.$$
Let also ${\cal P} = -i{\cal H}$ be the generalized right half-plane,
so it is given by the condition $x+x^*\geq 0$.
The action of the group $G\times G$ on $E$ does not preserve Herm, but the
action of the real Lie subgroup $\Delta = \{(g^{-1}, g^*)\in G\times G\}$ does.
As a real Lie group, $\Delta\simeq G$ but its embedding into $G\times G$ is 
not
holomorphic. Let {\bf d} be the (real) Lie algebra of $\Delta$. Thus {\bf d} 
is a real Lie
subalgebra in the complex Lie algebra ${\bf g}\oplus {\bf g}$ and ${\bf 
d}\otimes {\bf C}=
{\bf g}\oplus {\bf g}$.

For a subset $U\i G$ we denote by ${\cal M}_\alpha(U)$ the space of sections 
of ${\cal M}_\alpha$
over $U$.  Since ${\rm Herm}_+$ and $\cal P$ are, topologically, 
cells containing $1$,
we have identifications
$${\cal M}_\alpha({\cal P}) \simeq {\cal M}_\alpha({\rm Herm}_+) \simeq {\cal 
M}_{\alpha, 1}.
\leqno (4.2.1)$$
The power function $t^z, t\in {\bf R}_+$, $z\in {\bf C}$, will
be always normalized by putting $t^z = e^{z\log(t)}$ with the
standard choice of branch of logarithm on ${\bf R}_+$. We will
understand the expression $|x|^z$ (power of the determinant)
for $x\in {\rm Herm}_+$ accordingly. 
This defines an identification ${\cal M}_{\alpha, 1}\simeq M_\beta$,
where $\beta$ is related to $\alpha$ as in (4.1.2).
Similarly, we define an identification
$${\cal M}_\alpha({\cal H}) \simeq {\cal M}_\alpha({\rm Herm}_+) \simeq {\cal 
M}_{\alpha, 1}
\simeq M_\beta \leqno (4.2.2)$$
via the inclusion of ${\rm Herm}_+$ into the closure of $\cal H$. 

Let $\alpha\in {\bf I}$ be such that $\Re (\alpha_n) > 0$. Define the linear 
operator (Fourier
integral)
$$F_\alpha: {\cal M}_\alpha ({\rm Herm}_+)\to {\cal O}({\cal H}), \quad
(F_\alpha f)(y) = \int_{{\rm Herm}_+} f(x) e^{i{\rm tr}(xy)} dx,\leqno 
(4.2.3)$$
the integral being absolutely convergent. By identifying $\cal H$ with $\cal 
P$ via
multiplication by $i$, we can write the same map in the form of a Laplace
integral:
$$\Lambda_\alpha: {\cal M}_\alpha({\rm Herm}_+) \to {\cal O}({\cal P}),
\quad (\Lambda_\alpha f)(y) =  \int_{{\rm Herm}_+} f(x) e^{-{\rm tr}(xy)} 
dx,\leqno (4.2.4)$$
Recall that the Lie algebra ${\bf g}\oplus {\bf g}$ acts on ${\cal M}_\alpha 
({\rm Herm}_+)$.
We will also use the following {\it twisted action} of ${\bf g}\oplus {\bf g}$
on ${\cal O}({\cal P})$: 
$$(\gamma_1, \gamma_2) *f = n\bigl( {\rm tr}(\gamma_2) - {\rm 
tr}(\gamma_1)\bigr) \cdot f
+ (\gamma_2, \gamma_1)f,$$
where the last summand on the right is the standard action. 

\proclaim (4.2.6) Proposition. The map $\Lambda_\alpha$ is ${\bf g}\oplus {\bf 
g}$-equivariant,
if we consider the standard action on ${\cal M}_\alpha({\rm Herm}_+)$ and the 
twisted action on
${\cal O}({\cal P})$. 

\noindent {\sl Proof:} Since $\Lambda_\alpha$ is {\bf C}-linear, it is enough 
to prove
 {\bf d}-equivariance. This follows via a linear change of variables on ${\rm 
Herm}_+$.

\proclaim (4.2.7) Theorem. The map $\Lambda_\alpha$ defines an isomorphism of 
${\cal M}_\alpha
({\rm Herm}_+)$ with
\hfill\break
 ${\cal M}_{-n-\alpha_n, ..., -n-\alpha_1}({\cal P}) \i {\cal O}({\cal P})$.

\noindent {\sl Proof:} The space ${\cal M}_\alpha ({\rm Herm}_+)$ being an 
irreducible
${\bf g}\oplus {\bf g}$-module, it is enough to find the effect of
$\Lambda_\alpha$ on a generator of it, for example on the function 
$s_\alpha(x) = {\rm tr}(x^\alpha)$.
This has been done in [23], Th. 5.9, which, after being translated to our 
notation, reads:
$$ (\Lambda_\alpha s_\alpha)(y) = \int_{{\rm Herm}_+} e^{-{\rm tr}(xy)} 
s_\alpha(x) dx = $$
$$=\pi^{n(n-1)\over 2} \prod_{j=1}^n \Gamma(\alpha_j+n-j+1)
s_{-n-\alpha_n, ..., -n-\alpha_1}(y).\leqno (4.2.8)$$
(See also [32], formula (12.6) for an evaluation of an equivalent
integral.)
So the generator is mapped into a nonzero scalar multiple of a generator and 
our statement follows.

\vskip .3cm

\noindent {\bf (4.3) The matrix gamma function.}  For $\alpha\in {\bf I}$ set
$$\Gamma_n(\alpha) = \prod_{j=1}^n \Gamma(\alpha_j+n-j)\leqno (4.3.1)$$
and call this meromorphic function on {\bf I} the 
{\it matrix gamma function}. For $s\in{\bf C}$ we will write
$\Gamma_n(\alpha+s)$ for $\Gamma_n(\alpha_1 +s, ..., \alpha_n+s)$.
The map $\Lambda_\alpha$, as it follows from (4.2.8)
and equivariance, is essentially scalar:
$$\Lambda_\alpha[(a, x^\alpha)] = \pi^{n(n-1)\over 2} \Gamma_n(\alpha+1)
(a^t, x^{-n-\alpha_n, ..., -n-\alpha_1}).
\leqno (4.3.2)$$
Here $a\in {\rm End}(V_\alpha)$ and $a^t\in {\rm End}(V_\alpha^*)$ is its 
transpose.
We can also express this as follows:
$$\Lambda_\alpha(x^\alpha) = \pi^{n(n-1)\over 2} \Gamma_n(\alpha+1)
(x^{-n-\alpha_n, ..., -n-\alpha_1})^t,
\leqno (4.3.3)$$
where we regard $x^\alpha$ as a $V_\alpha^*\otimes V_\alpha$-valued function, 
while
$ x^{-n-\alpha_n, ..., -n-\alpha_1}$ is a $V_\alpha\otimes V_\alpha^*$-valued 
function.

\vskip .3cm

\noindent {\bf (4.4) Distributional version.}
Since the Fourier transform is involutive, a different setting will lead to 
the appearance of
$\Gamma_n(\alpha+1)^{-1}$ and not of $\Gamma_n(\alpha+1)$ as a proportionality 
factor in the
Fourier transform of $x^\alpha$. To overview both settings, it is convenient to
work with distributions on Herm, directly generalizing [20], \S II.2.

First, define a  ${\rm End}(V_\alpha)$-valued distribution $x^\alpha_+$
on Herm given by
$$x^\alpha_+ = \cases{x^\alpha, & if $x\in {\rm Herm_+}$;\cr
0&otherwise.\cr}\leqno (4.4.1)$$
Second, since Herm is the Shilov boundary of the Siegel upperhalf plane
 $\cal H$, any analytic
function $f$ in $\cal H$ has a boundary value on Herm, which
is a hyperfunction  denoted by
$$\bar f(x) = f(x+i0) = \lim_{\tau\in {\rm Herm_+}\atop \tau\to 0} 
f(x+i\tau).\leqno (4.4.2)$$
The identification ${\cal M}_\alpha({\cal H}) \simeq {\cal M}_{\alpha, 1}$ from
 (4.2.2), gives an ${\rm End}(V_\alpha)$-valued analytic function $x^\alpha$ 
in $\cal H$.
Its boundary value $(x+i0)^\alpha$ on Herm is in fact a distribution. 

Now, the Fourier transform of distributions on Herm is induced by the transform
 of test functions
given by
$$F[\phi](y) = \int_{\rm Herm} \phi(x) e^{i{\rm tr}(xy)} dx. \leqno (4.4.3)$$
It is also denoted by $F$ and satisfies 
$$FF[f(x)] = (2\pi)^{n^2} f(-x). \leqno (4.4.4)$$
Now, (4.3.3) implies, similarly to [20], 
the following equality of matrix-valued
distributions on Herm:
$$F[x_+^\alpha] = e^{i\pi(|\alpha|+n)\over 2} \pi^{n(n-1)\over 2} 
\Gamma_n(\alpha+1) (y^t+i0)^{(-n-\alpha_n, ..., -n-\alpha_1)},
\leqno (4.4.5)$$
and therefore
$$F[(x+i0)^\alpha] = 2^{n^2} \pi^{n(n+1)\over 2} { \biggl( y_-^{(-n-\alpha_n, 
...,
-n-\alpha_1)}\biggr)
^t\over \Gamma_n(-n-\alpha_n+1, ..., -n-\alpha_1+1)}.\leqno (4.4.6)$$

\vskip .2cm

\noindent {\bf (4.5) Divided powers.} Let $\alpha\in {\bf I}$. We denote
$$x^{((\alpha))} = {x^\alpha\over \Gamma_n(\alpha+1)}, \leqno (4.5.1)$$
regarding this as a matrix-valued function in the right half-plane $\cal P$.
Here we use the normalization of $x^\alpha$ in $\cal P$ described in (4.2.1).
The function $x^{((\alpha))}$ will be called the divided power. 

It follows from Proposition 1.3.3 that the  derivatives of $x^{((\alpha))}$
 have a particularly simple
form. In the next proposition we conserve the notation of (1.3).

\proclaim (4.5.2) Proposition. For $a\in {\rm End}(V^*), b_\alpha\in {\rm 
End}(V_\alpha)$
 we have
$$(a,\partial)\cdot (b_\alpha, x^{((\alpha))}) = 
\sum_{\beta=\alpha -e_i} \biggl( (p_{(1), \beta}^\alpha)^t
(a\otimes b_\alpha),\, x^{((\beta))}\biggr).$$

\vskip .2cm

\noindent {\bf (4.6) Example: exponential series.} If  $\alpha\in\Irr, 
\alpha\geq 0$ is
 a positive
Young diagram with $|\alpha|=m$, then let $W_\alpha$ be the corresponding 
representation
of the symmetric group $S_m$. Let $w_\alpha=\dim(W_\alpha)$. By the Weyl 
reciprocity theorem,
$$V^{\otimes m} \simeq \bigoplus_{|\alpha|=m} W_\alpha\otimes V_\alpha \leqno 
(4.6.1)$$
as a $G$-module, which implies that
$${\rm tr}(x)^m = \sum_{|\alpha|=m} w_\alpha \cdot s_\alpha(x).\leqno (4.6.2)$$
For any Young diagram $\alpha = (\alpha_1 \geq ... \geq \alpha_m)$
with $m$ cells the number $w_\alpha$ is given by the well known formula 
(see [11], formula (4.11)):
$$w_\alpha = { m!\over \prod_{j=1}^m (\alpha_j+m-j)!} \prod_{1\leq i<j\leq m} 
(\alpha_i-\alpha_j+j-i).
\leqno (4.6.3)$$
Note that as long as $\alpha$ has no more than $n$ non-zero parts
 (which is true for $\alpha$ entering
into (4.6.2)),
$${ \prod_{1\leq i<j\leq n}(\alpha_i-\alpha_j+j-i)\over \prod_{j=1}^n 
(\alpha_j+n-j)!}
= 
{ \prod_{1\leq i<j\leq m}(\alpha_i-\alpha_j+j-i)\over \prod_{j=1}^m 
(\alpha_j+m-j)!}.$$
Therefore, using 
the formula (1.1.7) for $d(\alpha)$ we get:
$${w_\alpha s_\alpha(x)\over |\alpha|!} = c_n {d(\alpha) 
s_\alpha(x)\over\Gamma_n(\alpha+1)},
 \quad c_n = \prod_{j=1}^{n-1} j!\leqno (4.6.4)$$
or, equivalently,
$$e^{{\rm tr}(x)} = 
\sum_{m=0}^\infty {{\rm tr}(x)^m\over m!} = c_n  \sum_{\alpha\in {\rm Irr}} 
d(\alpha)
 {s_\alpha(x)\over \Gamma_n(\alpha+1)}.\leqno (4.6.5)$$
Here we extended the sum over the entire set Irr, since 
$1/\Gamma_n(\alpha+1)=0$
 unless $\alpha\geq 0$. 
The equality (4.6.4) can be seen as the termwise Fourier transform of the 
formula
$$\delta_1(x) = \sum_{\alpha\in\Irr} d(\alpha) s_\alpha(x),\leqno (4.6.6)$$
discussed in (3.4.1). Note that (4.6.6) is an equality of distributions
on the unitary group, so it is not possible to formally deduce
(4.6.5) from it, if in the Fourier transform we integrate
over the space of Hermitian matrices.

\vskip .3cm

\noindent {\bf (4.7) Contour Fourier transform of single-valued
monomials.} When considering single-valued functions on $G=GL_n$,
we can take the contour in the Fourier integrals to be the
unitary group $U_n$. Thus, we define the {\it contour Fourier transform}
to be the map
$${\rm FC}: {\bf C}[G]\to {\bf C}[G], \quad {\rm FC}[f] (y)
=
\int_{x\in U_n} f(x) e^{{\rm tr}(xy)} dx.\leqno (4.7.1)$$
Here the integral can be, of course, calculated purely
algebraically, in terms of coefficients of the expansion of $f$
into monomials. The analog of the formula (4.4.6) in this algebraic
situation is as follows:
$${\rm FC}[x^{(-n-\alpha_n, ..., -n-\alpha_n)}] = C\cdot {y^\alpha\over
\Gamma_n(\alpha+1)}, \leqno (4.7.2)$$
where $C$ is the same as in (3.2.1). The formula (4.7.2) follows at once from
(3.2.2) and from the series expansion (4.6.5) of $e^{{\rm tr}(x)}$.

\vskip 2cm 

\centerline {\bf \S 5. The $A$-hypergeometric system.}

\vskip 1cm

\noindent {\bf (5.1) The system and its Fourier transform.} We consider
 a pair of reductive groups
$$H\i G = \prod_{\omega\in A} GL(V_\omega), \quad A\i \Irr(H),$$
satisfying the homogeneity property (2.1.3). We define the space
$M_A^* = \bigoplus_{\omega\in A} {\rm End}(V_\omega)$ and varieties
$Y_A\i M_A^*, X_A\i P(M_A^*)$ as in \S 2. Thus $Y_A$ is 
$H$-biinvariant and contains $H$ as an open orbit. Let $I_A\i S^\bullet
(M_A)$ be the homogeneous ideal of $Y_A$. 

Let {\bf h} be the Lie algebra of $H$, and $\chi: {\bf h}\to {\bf C}$ be
a character. Let us identify $S^\bullet(M_A)$ with the ring of differential
operators on $M_A$ with constant coefficients. More precisely, for 
$f\in S^\bullet(M_A)$ we denote by $P_f$ the corresponding differential
operator. 

By definition, the $A$-{\it hypergeometric system} corresponding to the 
character
$\chi$, is the following system of linear differential equations on a function
$\Phi\in {\cal O}(M_A)$:
$$\cases{ L_h\Phi = R_h\Phi  = \chi(h)\cdot \Phi, \quad h\in{\bf h},\cr
P_f\Phi = 0, \quad f\in I_A.\cr} \leqno (5.1.1)$$
Here $L_h$ and $R_h$ are the infinitesimal generators of the 
left and right $H$-actions on $M_A$.
Thus the first group of equations just expresses the quasihomogeneity of 
$\Phi$.
Holomorphic solutions of this system (defined over some open sets in $M_A$)
will be called $A$-hypergeometric functions. 
We denote by ${\cal H} = {\cal H}_{A,\chi}$ the sheaf of such functions.  

\vskip .1cm

For a finite-dimensional space $V$ let $\D(V)$ denote the algebra of 
differential
operators on $V$ with polynomial coefficients. There is a natural isomorphism
$\phi: \D(V)\to \D(V^*)$ called the (formal) Fourier transform. If $x_1, ..., 
x_N$ are
linear coordinates in $V$ and $y_1, ..., y_N$ are the dual coordinates in 
$V^*$, then
$$\phi(x_i) = -{\partial\over\partial y_i}, \,\, 
\phi\biggl({\partial\over\partial x_i}
\biggr) = y_i.$$
For $P\in \D(V)$ we will denote $\phi(P)$ simply by $\hat P$.

Now take $V=M_A$ and apply the Fourier transform to the system (5.1.1). In 
this way
we get the following system on a function
$\Psi$ on $M_A^*$;
$$\cases{ \hat L_h \Psi = \hat R_h\Psi =\chi(h)\cdot \Psi, \quad h\in {\bf 
h}\cr
f\cdot \Psi = 0, \quad f\in I_A.\cr} \leqno (5.1.2)$$
Note that  $\hat L_h$, resp. $\hat R_h$ is just the infinitesimal generator of 
the
left, resp. right  $H$-action on $M_A^*$
dual to that on $M_A$. The second group of equations means that $\Psi$
is supported on $Y_A$, so solutions exist only among distributions.

\proclaim (5.1.3) Theorem. (a) The system (5.1.1) is holonomic, with regular 
singularities.
In particular, there is an open subset $M_A^{gen}$ in $M_A$ (``generic 
stratum") such
that the restriction of the sheaf $\cal H$ to  $M_A^{gen}$ is a locally 
constant sheaf on
$M_A$ of finite rank.\hfill\break
(b) More precisely, $M_A^{gen}$ is obtained from $M_A$ by deleting the
varieties $\nabla_{A, \Gamma}$, see (2.8), for all proper faces 
$\Gamma$ of the polytope $Q_A$. 
\hfill\break
(c)   Supose that the variety $Y_A$ is Cohen-Macaulay.
Then the rank of the sheaf $\cal H$ at any point of $M_A^{gen}$ is less or 
equal to the degree
of $X_A$ (which is given by the Kazarnovskii-Brion theorem 2.7.1). 

\noindent {\bf (5.1.4) Remark.} In the toric case, the Cohen-Macaulay
property implies that the generic rank of $\cal H$ is equal to
${\rm deg}(X_A)$, see [15]. Without this property, the generic rank
of $\cal H$ can be greater, as shown by B. Sturmfels and N. Takayama [47].

\vskip .3cm

\noindent {\bf (5.2) The hypergeometric $\cal D$-module.} The proof of (5.1.3) 
is similar to
that in [15] and is based on the analysis of the $\cal D$-module corresponding 
to the system
(see [2] [28] for  general background on $\D$-modules). 
If $Z$ is a smooth algebraic variety, we denote by ${\cal D}_Z$ the sheaf of 
regular differential
operators on $Z$ and by ${\cal D}(Z)$ the ring of its global sections. 
For a coherent $\D_Z$-module $\cal N$ we denote by $SS({\cal N})\i T^*Z$ its 
characteristic
variety and by ${\bf SS}({\cal N})$ the characteristic cycle.

We denote for short $\D = \D_{M_A}$ and $\hat\D = \D_{M_A^*}$
Let ${\cal M} = {\cal M}_{A,\chi}$ be the $\D$-module corresponding to the 
system (5.1.1):
$${\cal M}_{A,\chi}= \D\biggl/ \sum_{f\in I_A} \D \cdot P_f + \sum_{h\in{\bf 
h}} \D\cdot
(L_h-\chi(h)) + \sum_{h\in{\bf h}} \D (R_h-\chi(h)).\leqno (5.2.1)$$
Thus the sheaf of $A$-hypergeometric functions can be written as
$${\cal H}_{A,\chi} = \underline{\rm Hom}_\D ({\cal M}_{A,\chi}, {\cal 
O}_{M_A}).
\leqno (5.2.2)$$

Let also $\hat {\cal M}$ be the $\hat\D$-module corresponding to the system 
(5.1.2).
It is a particular case of $\D$-modules described in the following theorem,
see [2], Ch. VII, Th. 12.11 or [28], \S 5.

\proclaim (5.2.3) Theorem. Let $K$ be an algebraic group over {\bf C} 
acting on
a smooth algebraic variety $Z$. Let $W\i Z$ be a closed
subvariety which is the union of  finitely many $K$-orbits and let ${\cal 
I}_W\i {\cal O}_Z$
be the sheaf of functions vanishing on $W$. Let ${\bf k}$ be the Lie algebra 
of $K$
and $\beta: {\bf k}\to {\bf C}$ be a character. Then the $\D_Z$-module
$$\D_Z\biggl/ \D_Z\cdot {\cal I}_W + \sum_{k\in{\bf k}} \D_Z(L_k-\beta(k))$$
is holonomic, with regular singularities. Its characteristic variety is 
contained in the
union of the closures   of the conormal bundles $T^*_O(Z)$, where $O$ runs 
over the $K$-orbits in $W$.

The module $\hat{\cal M}$ corresponds to $K=H\times H$, $Z=M_A^*$, $W=Y_A$.  
The proof of Theorem 5.1.3 (a-b) is now done in the same way as in [15]: the 
module
$\hat {\cal M}$ being monodromic in the sense of [6], its Fourier transform, 
i.e., $\cal M$, is again
holonomic regular by Corollary 7.25 of [6], which gives (a). To see (b), let
$\bar\phi: T^*M_A\to T^*M_A^*$ be the natural identification of the contangent 
spaces.
The result  of [6] just cited
implies also that $SS({\cal M}) = \bar\phi^{-1}(SS(\hat{\cal M}))$ which gives 
(b).

\vskip .1cm

We now prove part (c). We start with the following lemma.

\proclaim (5.2.4) Lemma. Let $Y\i {\bf C}^N$ be  a conic (${\bf 
C}^*$-invariant)
closed algebraic variety, whose projectivization $P(Y)\i P^{N-1}$ has degree
$d$, and let $E\i {\bf C}^N$ be a linear subspace such that
$\dim(E)+\dim(Y)=N$ and $Y\cap E=\{0\}$. If $Y$ is Cohen-Macaulay,
then
$$\dim_{\bf C} \biggl( {\bf C}[E]\otimes_{{\bf C}[x_1, ..., x_N]} {\bf C}[Y]
\biggr) = d.$$

\noindent{\sl Proof:} Let $P^N$ be the projective space containing ${\bf C}^N$,
and $\bar E, \bar Y$ be the closures of $E,Y$ in $P^N$,
so that $\bar E$ is a projective subspace and ${\rm deg}(\bar Y)=d$. Then
$\bar E\cap \bar Y= E\cap Y$ consist of one point, namely $0\in {\bf C}^N$.
Therefore
$$\sum_{i} (-1)^i \dim_{\bf C} {\rm Tor}_i^{{\bf C}[x_1, ..., x_n]} ({\bf C}
[E], {\bf C}[Y]) = d,$$
as this is the contribution of the only intersection point of $\bar E$
and $\bar Y$ into their global intersection index on  $P^N$, which is
$d={\rm deg}(\bar Y)$, see [43]. On the other hand, let $l_1, ..., l_m$
be independent linear equations of $E$. As $\dim(E\cap Y)=0 = \dim(Y)-m$,
we conclude that $l_1, ..., l_m$ form a regular sequence in ${\bf C}[Y]$,
see [38], Th. 17.4. The Koszul complex associated with this sequence
is precisely the complex calculating the Tor's above. Since
its homology in degrees other than 0 vanishes,
 we conclude that ${\rm Tor}_i=0$
for $i\neq 0$ and the lemma follows. 

\vskip .1cm

Let now $h\in{\bf h}$. The vector field $L_h$ on $M_A$ can be regarded as
a family of linear functionals $L_h(a): M_A^*\to {\bf C}$ depending 
on $a\in M_A$. These functionals form the highest symbol of the
equation $L_h\Phi = \chi(h)\cdot\Phi$ from (5.1.1).
For given $a\in M_A$ let $\lambda_h(a)\i M_A^*$ be the kernel of
$L_h(a)$. Similarly, let $R_h(a): M_A^*\to {\bf C}$ be the
linear functional corresponding to $R_h$, and $\rho_h(a)$ be its kernel.
Now, the following fact is clear from the definitions.

\proclaim (5.2.5) Lemma. Let $a\in M_A, h\in{\bf h}$. A point $b\in M_A^*$
lies in $\lambda_h(a)$, resp. in $\rho_h(a)$, if and only if the vector
$${d\over dt}\biggr| _{t=0} e^{th}\cdot b, \quad {\rm resp.}\quad 
{d\over dt}\biggr| _{t=0} b\cdot e^{th}\in M_A$$
is orthogonal to $a$. 

Let
$$\lambda(a) = \bigcap_{h\in {\bf h}} \lambda_h(a), \quad \rho(a) = 
\bigcap_{h\in{\bf h}} \rho_h(a)$$
be the zero loci over $a$ of the highest symbols of the first two
groups of equations in (5.1.1).  

 \proclaim (5.2.6) Corollary. A point $b\in M_A^*$ lies in $\lambda(a)$, 
resp. $\rho(a)$, resp. $\lambda(a)\cap\rho(a)$ if and only if
the hyperplane ${\rm Ker}(a)\i M_A^*$ is tangent to the left
$H$-orbit, resp. right $H$-orbit, resp. $H\times H$-orbit of $b$.

In particular, for $x_0\in H$ the point $b=(x_0^\omega)_{\omega\in A}$
of $\rho_A(H)\i M_A^*$ lies in $\lambda(a)$ or $\rho(a)$ if and
only if the polynomial $f_a(x) = \sum (a_\omega, x^\omega)$ on $H$
has a critical point at $x=x_0$ (the critical value is automatically
0 by homogeneity).

\proclaim (5.2.7) Corollary.  If $a\in M_A^{gen}$, then $\lambda(a)\cap\rho(a)
\cap Y_A=\{0\}$. 

Now, noticing that ${\rm codim} \,\,\lambda(a) = {\rm codim}\,\,\rho(a) =
\dim Y_A$, we can find a linear subspace $E\i M_A^*$ of codimension equal
to $\dim(Y_A)$, containing $\lambda(a)\cap\rho(a)$
and such that $E\cap Y_A=\{0\}$. Therefore, by Lemma 5.2.4, the quotient of 
${\bf C}[Y_A]$ by the ideal generated by the linear equations of $E$, has 
dimension
equal to ${\rm deg}(X_A)$. Thus the quotient by the ideal generated by
the linear equations of $\lambda(a)\cap\rho(a)$ has dimension less
or equal to ${\rm deg}(X_A)$. Finally, this latter quotient is just
obtained by taking the highest symbols of the equations in (5.1.1),
so the quotient by the full characteristic ideal has the rank less or
equal to that. This proves part (c) of Theorem 5.1.3. 
\vskip .3cm

\noindent {\bf (5.3) Fourier integrals. }  We keep the above assumptions on 
$H, A, \chi$.
The system of differential equations
$$L_h(u) = \chi(h)\cdot u, \quad h\in {\bf h}, \leqno (5.3.1)$$
together with the initial condition $u(1)=1$, defines a multivalued 
``monomial" function
on $H$ which we will denote $x^\sigma$, $\sigma =\sigma(\chi)$. Denote by 
${\cal L}={\cal L}_\chi$
the 1-dimensional local system on $H$ spanned by all the branches $x^\sigma$. 
Let $C_\bullet(H, {\cal L})$ be the singular chain complex of $H$ with 
coefficients in $\cal L$,
see, e.g., [16]. An element of $C_m(H, {\cal L})$ is thus a finite linear 
combination of
pairs $(\gamma, \phi)$ where $\gamma: \Delta^m\to H$ is a singular $m$-simplex 
and $\phi$ is a section of
$\gamma^*{\cal L}$ on $\Delta^m$. The differential is defined by the usual 
rule. Let
$C^{lf}_\bullet(H,{\cal L})$ be the bigger complex consisting of locally finite
combinations as above. As well known, the homology of these complexes is 
expressed via sheaf
cohomology (with or without compact support):
$$H_m(H, {\cal L}) = H^m(H, {\cal L}^*)^*, \quad H_m^{lf}(H,{\cal L}) = 
H^m_c(H, {\cal L}^*)^*.$$
We regard $H$ as embedded into $M_A^*$. Let $K\i M_A^*$ be any strictly convex 
closed cone.
Denote by $C_\bullet^{(K)}(H,{\cal L})$ the subcomplex in $C_\bullet^{lf}(H, 
{\cal L})$ formed
by such locally finite chains which are contained in $K$ except for finitely 
many simplices
and which are semi-algebraic at the infinity of $M_A^*\supset H$ in the
sense of [40]. 
Let also $C_\bullet^{[K]}(H,{\cal L})$ be the bigger  subcomplex consisting of 
such locally
finite chains that all their simplices that do not belong to $K$,
 lie in a bounded subset of $M_A^*$ (and which satisfy the same
 condition at the infinity of $M_A^*$ as before). 
We will denote the homology of these complexes by $H_m^{(K)}(H, {\cal L})$ 
etc.

Let $N=\dim(H)$. We define the complexes of sheaves ${\cal F}, {\cal F}^!$ on 
$M_A$ by their stalks:
$${\cal F}_a = \lim_{K\i \{ \Im(b,a) > 0\}} H^{(K)}_N (H, {\cal L}), \quad
{\cal F}^!_a = \lim_{K\i \{ \Im(b,a) > 0\}} H^{[K]}_N (H, {\cal L}), \quad 
a=(a_\omega)\in M_A.
\leqno (5.3.2)$$
Here $\Im$ means the imaginary part of a complex number. 
Define the morphism of sheaves (Fourier integral)
$${\cal I}: {\cal F}\to {\cal O}_{M_A}, \quad \zeta=\sum_\nu (\gamma_\nu, 
\phi_\nu) \mapsto
\sum_\nu \int_{\gamma_\nu}\phi_\nu (b)
 \exp\biggl(i\sum {\rm tr}(a_\omega b_\omega)\biggr) db.\leqno (5.3.3)$$
Here the sum (i.e., the improper integral) always converges because of the
 fast decay of the exponent
in any cone $K$ in (5.3.2). Note that in (5.3.3) we regarded $H$ as a 
subvariety in $M_A^*$
given by the parametric equations $b_\omega = x^\omega, x\in H$. Thus, 
regarding the
linear function $\sum {\rm tr} (a_\omega b_\omega)$ as a polynomial
$f(x) = \sum (a_\omega, x^\omega)$ on $H$, we can rewrite the map $\cal I$ as 
follows:
$${\cal I}\biggl(\sum_\nu (\tau_\nu, \phi_\nu)\biggr)(f) = \sum_\nu 
\int_{\tau_\nu}
\phi_\nu e^{if(x)}dx.\leqno (5.3.4)$$

\proclaim (5.3.5) Proposition. The image under $\cal I$ of any section of 
$\cal F$ satisfies
the $A$-hypergeometric system. 

\noindent{\sl Proof:} Straightforward.

\vskip .2cm

\noindent {\bf (5.3.6) Remark.} 
The above proposition is just a manifestation of the general principle that the
Fourier transform of a solution of a linear differential system, whenever it 
makes
sense, is a solution of the (formally) Fourier  transformed system.
More precisely, since $H=\rho_A(H)$ is an open orbit in $Y_A$, the space of 
solutions
of (5.1.2) near the generic point of $Y_A$  is 1-dimensional and generated
by the following ``holomorphic distribution":
$$(\rho_A)_*(x^\sigma) = \delta_H(b) \cdot \rho_A^{-1}(b)^\sigma, \quad b\in 
G\i M_A^*.
\leqno (5.3.7)$$
This is just the delta function along $H\i G$ multiplied by a power function.
The sheaf $\cal F$ is formed by all possible contours which can be used to 
define the
Fourier transform of this delta function. A more suggestive and classical way 
to rewrite
a function from the image of $\cal I$ is:
$$\Phi(f) = \int_{\zeta \i H} x^\sigma e^{if(x)}dx = \int_{\zeta\i Y_A\i 
M_A^*} \exp\biggl(
i\sum {\rm tr}(a_\omega\cdot b_\omega)\biggr) \delta_H(b) 
\rho_A^{-1}(b)^\sigma db.\leqno (5.3.8)$$
Here $\zeta$ is an ``appropriate contour of integration" (the exact meaning of 
this is
provided by the sheaf $\cal F$). For the toric case Fourier integral
representations of hypergeometric functions were studied in [12-13]
where the case of an arbitrary locally compact ground field $k$ is considered.
For $k={\bf R}$ (and $H$ a torus) 
the framework of [12-13] gives the restrictions
of $A$-hypergeometric functions to  real values of variables while
for $k={\bf C}$ it is slightly different from the one presented here
in that we allow more complicated domains of integration.

\vskip .2cm

\noindent {\bf (5.4) Euler integrals.} 
 Suppose now that we have a pair $(H_0, A_0)$ satisfying 
the inhomogeneity condition (2.1.4) and let $(H, A)$ be the group
and the set of representations obtained in (2.1), so that they satisfy the
homogeneity condition and the $A$-hypergeometric system is defined.
We have an exact sequence
$$1\to H_0\cap {\bf C}^*\to H_0\times {\bf C}^*\buildrel p\over\to  H
\to 1, \quad  H\i G = \prod_{\omega\in A_0} GL(V_\omega)
\leqno (5.4.1)$$
with $p$ being the multiplication map. Recall also that we have the 
identification
$A_0\to A$, $\omega\mapsto \bar\omega$ with $V_{\bar\omega}=
V_{\omega}$ as a vector space, and thus $M_{ A_0} = M_A$. 
We  want to represent the solutions of the $A$-hypergeometric
system (which are functions on $M_{ A}=M_{A_0}$) as certain integrals
of $f\in M_{A_0}$ over $H_0$. 

More precisely, let $\chi$ be a character of ${\bf h} = {\rm Lie}( H)$. 
Since the map $p$ in (5.4.1) is finite, we have an identification
${\bf h} = {\bf h}_0\oplus {\bf C}$. The character $\chi$ is thus a pair
$(\chi_0, \tau)$, where $\chi_0$ is a character of {\bf h} and $\tau\in {\bf 
C}$.
Let  $x^{\sigma_0}$, $x\in H_0$, be the multivalued group character
corresponding to $\chi_0$, similarly to (5.3.1). 

For $f\in M_A$ we denote by $U_f$ the open set $\{f\neq 0\}\i H$. Let
${\cal N}_\chi$ be the 1-dimensional local system on $U_f$ spanned by all the 
branches of the
 multivalued function $x^{\sigma_0} f(x)^\tau$. 
Let $N_0=\dim (H_0)$. Consider the constructible sheaf ${\cal E}={\cal 
E}(\chi)$
on $M_A$ defined by its stalks as follows:
$${\cal E}_f = H_{N_0}(U_f, {\cal N}_\chi). \leqno (5.4.2)$$
We have the morphism of sheaves on $M_A$
$${\cal J}: {\cal E}\to {\cal O}_{M_A}, \quad \zeta_0  = \sum_\nu (\gamma_\nu, 
\phi_\nu)\mapsto
\sum_\nu \int_{\gamma_\nu} \phi_\nu(x) dx. \leqno (5.4.3)$$
Images of sections of $\cal E$ under this map will be called {\it Euler 
integrals}.

\proclaim (5.4.4) Proposition. Every Euler integral satisfies the 
$A$-hypergeometric
system. 

The proof is straightforward.

\vskip .1cm

\noindent {\bf (5.4.5) Remark.} It is more suggestive to write an
 Euler integral
in the more classical form:
$$J(f) = \int_{\zeta_0} x^{\sigma_0} f(x)^\tau dx,\leqno (5.4.6)$$
where $\zeta_0\i H_0$ is ``an appropriate contour of integration". 
On the heuristic level, the Euler integral over $\zeta_0\i H_0$
 is obtained from a Fourier integral over the cycle $\zeta = \zeta_0 \times
(i{\bf R}_+)$ in $H_0\times {\bf C}^*$ by integrating away the ${\bf 
R}_+$-variable:
$$\int_{\zeta} e^{i\lambda f(x)}\lambda^\tau x^{\sigma_0}
dx d\lambda = \int_{\zeta_0} \biggl( \int_{{\bf R}_+} e^{-\lambda f(x)} 
\lambda^\tau d\lambda
\biggr)
x^{\sigma_0} dx 
=\Gamma(\tau +1) \int_{\zeta_0} f(x)^{-\tau - 1} x^{\sigma_0} dx.$$
On the other hand, the first integral is equal to
$$ \int_{p(\zeta)\i  H} \exp(i{ \bar f}(\bar x)) \bar x^\sigma d\bar x,$$
where $\bar f\in M_A$ corresponds to $f\in M_{A_0}$ under the canonical 
identification.

\vskip .3cm

\noindent {\bf (5.5) Nonresonance.} 
As in (2.4), choose a maximal torus $T\i H$, let $\Lambda$ be its lattice
of characters and $\Lambda_{\bf R}=\Lambda\otimes {\bf R}$.
We assume other notations of (2.4) as well. Let $C_A\i\Lambda_{\bf R}$ be
the convex cone with apex $0$ and base $Q_A$. Thus nonempty faces of $C_A$
are in bijection with faces of $Q_A$. 

Let $\Lambda^0\i\Lambda$ be the lattice of characters of the full group
$H$, and $\Lambda^0_{\bf C}=\Lambda^0\otimes {\bf C}$. Thus the character
$\chi$ defining the $A$-hypergeometric system, lies in $\Lambda^0_{\bf C}$. 
For a face $\Gamma\i C_A$ we denote ${\rm Lin}_{\bf C}(\Gamma)\i\Lambda_{\bf 
C}$
its {\bf C}-linear span. We will say that $\chi\in\Lambda^0_{\bf C}$
is nonresonant, if  for any face $\Gamma\i C_A$ of codimension 1 we have
$\chi\notin \Lambda + {\rm Lin}_{\bf C}(\Gamma)$. 

Note that unlike the toric case [16], it may be no longer the case that a 
generic
$\chi\in\Lambda^0_{\bf C}$ is nonresonant. This will be true, however,
if the set $A$ is chosen generic enough (so that $\Lambda^0$ does not
lie in any ${\rm Lin}_{\bf C}(\Gamma), \Gamma\i C_A$).

The meaning of the nonresonance property is as follows.

\proclaim (5.5.1) Proposition. If $\chi$ is nonresonant, then the natural
morphism  $j_!{\cal L}_\chi \to Rj_* {\cal L}_\chi$ is an isomorphism
in the derived category of sheaves on $M_A^*$. 

\noindent {\sl Proof:} The closure of $H$ in $M_A^*$ is $Y_A$, and
$H\times H$-orbits $Y(\Gamma)$ on $Y_A$ are in bijection with
nonempty faces $\Gamma\i C_A$. 
We have to prove the following statement: for each proper face $\Gamma$,
each point $y\in Y(\Gamma)$ and small neighborhood $U$ of $y$ in
$Y_A$ we have $H^i(U\cap H, {\cal L}_\chi) =0$ for all $i$.

Consider first the case ${\rm codim}(\Gamma)=1$. In this case
${\rm codim}(Y(\Gamma))=1$, since the open orbit $H\i Y_A$ is affine.
Applying the slice theorem ([3], n. 1.5) to the normalization of $Y_A$,
we find that for $y,U$ as above, $U\cap H$ is homotopy equivalent
to the disjoint union of several copies of the circle $S^1$. The condition
$\chi\notin \Lambda + {\rm Lin}_{\bf C}(\Gamma)$ means that the monodromy
of ${\cal L}_\chi$ on each of these circles is nontrivial, so all the 
cohomology
vanishes.

The case ${\rm codim}(\Gamma) >1$ follows from this, since for such $\Gamma$
and $U,y$ as above, $U\cap H$ is fibered into  unions of circles corresponding 
to
any $\Gamma'\supset \Gamma$, ${\rm codim}(\Gamma')=1$. 

\vskip .1cm

It seems that for Cohen-Macaulay $Y_A$, the nonresonance condition should
be sufficient to ensure that the systems of Euler and Fourier integrals are
complete, i.e., the maps ${\cal I}: {\cal F}\to {\cal H}$
and ${\cal J}: {\cal E}\to {\cal H}$ are isomorphisms of sheaves.
A possible approach to this (generalizing that of [16]) would be
to find the characteristic cycle of the perverse sheaf $j_!{\cal L}_\chi$.
This is a purely topological problem. For example, the multiplicity of
$T_0^*M_A^*$ in the characteristic cycle is just the topological Euler
characteristic of the hypersurface $\{f=\epsilon\}\i H$ for generic
$f\in M_A, \epsilon\in {\bf C}$. When $H$ is a torus, this Euler
characteristic is equal to ${\rm deg}(X_A)$, as follows from [34].
It is likely that a similar answer is possible for an arbitrary reductive
$H$.

\vskip .3cm

\noindent {\bf (5.6)  Matrix $\Gamma$-series.} Another approach to solving the
hypergeometric system is to expand the holomorphic distribution (5.3.7)
 into a formal
power series and perform the  termwise Fourier transform,
using the formulas of (4.4).  Since we will work with power
series on $G=\prod_{\omega \in A} GL(V_\omega)$, let us introduce some 
notation.

An element of $G$ will be written as $a=(a_\omega)$. 
Let $J_\omega \simeq ({\bf C}^*)^{d(\omega)}$ be a maximal torus in 
$GL(V_\omega)$,
and $\Xi_\omega \simeq {\bf Z}^{d(\omega)}$ be its lattice of characters.
An element of $\Xi_\omega$ will be written as $\alpha(\omega) = 
(\alpha(\omega)_1, ..., \alpha(\omega)_{d(\omega)})$.  Let us abbreviate
$\Irr_\omega :=\Irr(GL(V_\omega))$ and identify this set with the set of
dominant weights in $\Xi_\omega$. A representation of $G$ will be written
as $V_{[\alpha]} = \bigoplus_{\omega\in A} \Sigma^{\alpha(\omega)}(V_\omega)$.
Let $\Irr_\omega^+$ be the set of positive dominant weights, i.e., those
with all $\alpha(\omega)_i\geq 0$.

Set $J=\prod_{\omega\in A} J_\omega$, $\Xi = \prod \Xi_\omega$,
$\Irr = \prod \Irr_\omega$.  Let also $\Xi^0 = {\bf Z}^A\i \Xi$ be the lattice
of characters of $G$. An element of $\Xi^0$ will be written as
$s=(s_\omega)$ and the corresponding character is $a\mapsto \prod
\det(a_\omega)^{s_\omega}$. 

We denote by $\Xi_{\bf R} = \Xi\otimes {\bf R}$ etc. 
Let also $\Irr_{\omega, {\bf R}}^+$ be the convex hull of $\Irr_\omega^+$ in 
$\Xi_{\omega,
{\bf R}}$. 

Denote by ${\bf I}_\omega$ the set of $\alpha(\omega) =
(\alpha(\omega)_{1}, ..., \alpha(\omega)_{ d(\omega)})\in {\bf C}^{d(\omega)}$ 
with
$\alpha(\omega)_i -\alpha(\omega)_{i+1}\in {\bf Z}_+$ and then set ${\bf I} = 
\prod_{\omega\in A}
{\bf I}_\omega$. Recall that the hypergeometric system depends on a character 
$\chi$ of {\bf h}.
Let $d\rho_\omega$ be the representation of {\bf h} corresponding to the 
representation
$\rho_\omega$ of $H$. The set {\bf I}, on the other hand, parametrizes 
irreducible representations
of ${\bf g}={\rm Lie}(G)$. We introduce the space
$$L_\chi = \biggl\{ s=(s_\omega)\in \Xi^0_{\bf C} =
 {\bf C}^A\biggl| \,\, \sum_{\omega\in A}  s_\omega {\rm tr} \,\,d\rho_\omega
 = \chi.\biggr\}.
\leqno (5.6.1)$$
In terms of the multivalued function $x^\sigma$ this condition is expressed 
just as
$\prod |x^\omega|^{s_\omega} = x^\sigma$. 

For any $s\in L_\chi$ the formal series
$$\sum_{\alpha\in\Irr} d(\alpha)\cdot \bigl(I_\alpha(H), a^{\alpha+s}\bigr), 
 \quad a^{\alpha+s} = \bigotimes
a_\omega^{\alpha(\omega)}\cdot |a_\omega|^{s_\omega},\leqno (5.6.2)$$
represents the distribution (5.3.7)
{\it on the compact form} $G_c\i G$. However, each term of this series
 also gives
a distribution on the space $M_A^{\rm Herm} = \prod _{\omega\in A} {\rm 
Herm}(V_\omega)$,
where ${\rm Herm}(V_\omega)\i {\rm End}(V_\omega) = M_\omega$ is the
space of Hermitian operators, see (4.4). The corresponding series
of distributions on $M_A^{\rm Herm}$ does not converge
but we still can form its 
  term-by-term
Fourier transform, using (4.4.6):
$$\Phi_s(a) = \sum_{\alpha\in \Irr} \biggl( { d(\alpha) \cdot I_\alpha(H)\over 
\prod_\omega \Gamma_{d(\omega)}
(\alpha_\omega+s_\omega+1)}, \,\, a^{\alpha+s}\biggr), \quad d(\alpha) = 
\prod d(\alpha_\omega).
\leqno (5.6.3)$$
We get in this way a {\it formal} series $\Phi_s$ called
the matrix $\Gamma$-series. 
It has the following meaning. let $U\i G$ be an open simply connected
domain, $a^s$ a branch of the function $\prod \det(a_\omega)^{s_\omega}$ in 
$U$,
and ${\cal M}_s$ the 1-dimensional space of functions in $U$ spanned by $a^s$.
Then $\Phi_s(a) \in {\bf C}[[G]]\otimes_{\bf C} {\cal M}_s$, where
${\bf C}[[G]]$ is the space of formal series on $G$ defined in (3.3).
The action of differential operators on elements of ${\bf C}[[G]]\otimes {\cal 
M}_s$
is defined similarly to Proposition 3.3.1. 

\proclaim (5.6.4) Proposition. (a) If $s'\in L_0 \cap {\bf Z}^A$, then 
$\Phi_{s+s'}(a) =
\Phi_s(a)$ as formal series.\hfill\break 
(b) The series $\Phi_s(a)$  satisfies the 
$A$-hypergeometric
 system in the sense described above.

\noindent {\sl Proof:} Part (a) is obvious,
part (b) follows from the properties of divided powers (Proposition 4.5.2).

\vskip .3cm

\noindent {\bf (5.7) Convexity and convergence of matrix $\Gamma$-series.} 
Let $K(G,H)\i \Xi_{\bf R}$ be the convex hull of all the weights of all
those representations $V_{[\alpha]}=\bigotimes \Sigma^{\alpha(\omega)}
(V_\omega)$ of $G$ for which $V_{[\alpha]}^H\neq 0$. This is a convex
but not necessarily strictly convex, cone in $\Xi_{\bf R}$. For a subset
$B\i A$ we denote $p_B: \Xi_{\bf R}\to \prod_{\omega\in B} \Xi_{\omega, {\bf 
R}}$
the coordinate projection. 

\proclaim (5.7.1) Definition. A subset $B\i A$ is called a cobase, if
$$K(G,H)\cap p_B^{-1}\biggl(\prod_{\omega\in B} \Irr_{\omega, {\bf 
R}}^+\biggr)$$
is a strictly convex cone. 

\noindent {\bf (5.7.2) Example.} Let $H$ be a torus, so that $\Irr(H)=\Lambda$
is a lattice. A subset $B\i A$ is a cobase if and only if the complement
$A-B$ is an affinely independent subset, i.e., forms the set of vertices of
a simplex, see [15]. 

\proclaim (5.7.3) Theorem. If $B$ is a cobase and $s\in \Xi^0_{\bf C}$ is such 
that
$s_\omega\in {\bf Z}$ for $\omega\in B$, then the series $\Phi_s(a)$ converges
in some open domain in $G$. 

If $s$ satisfies the conditions of the theorem, we shall say that $\Phi_s(a)$
is adapted to $B$.

\vskip .2cm

\noindent {\sl Proof:} Because of the poles of  the gamma-functions
entering the factors $\Gamma_{d(\omega)}(\alpha(\omega)+s_\omega+1)$ in the
denominators, we find that  $\Phi_s(a)$ is the product of a scalar
monomial and a series $\Phi'_s(a)$ whose support  lies
in the cone described in Definition 5.7.1, which is, by assumption, strictly
convex.  Denote by $S$ the set of all irreducible representations $V$
of $G$ such that all the weights of $V$
lie in this cone. Then $\Phi'(a)$ belongs to the ring
$M[[S]]$  of formal power series from (3.3) and satisfies a holonomic
system with regular singularities. Thus it has a nonempty domain
of convergence, by  Theorem  3.3.2.

\vskip .3cm

\noindent {\bf (5.8) Terminating series.} The set $B=A$
is always a cobase. Unlike the torus case, it is not possible to guarantee
the existence of any other cobases for  general $(H,A)$  (see \S 6 for examples
where nontrivial cobases exist).  For $B=A$ the condition of  Theorem 5.7.2
is simply that all the $s_\omega$ are integers. A matrix $\Gamma$-series
with this property
will be called totally resonant.

Let us identify ${\bf Z}^A$ with ${\rm Hom}(G, {\bf C}^*)$, associating to
$s=(s_\omega)$ the character $|b|^s=\prod |b_\omega|^{s_\omega}$.
Let $${\bf Z}^A \buildrel r\over\to \Lambda^0={\rm Hom}(H, {\bf C}^*)
\buildrel q\over\to {\rm Hom}({\bf C}^*, {\bf C}^*)={\bf Z}
\leqno (5.8.1)$$
be the restriction maps. We think of $\Lambda^0$ as a lattice
 (inside $\Irr(H)$) and
for an element $\sigma$ of this lattice denote by $x^\sigma$ the corresponding
character of $H$. Thus, the following is obvious.

\proclaim (5.8.2) Proposition. (a) If $\chi$ is a character of ${\bf h}$
and $s\in L_\chi\cap {\bf Z}^A$, $\sigma = r(s)$, then $x^\sigma$ is the
same as the function defined by (5.3.1). \hfill\break
(b) For $s\in {\bf Z}^A$ we have $q(r(s)) = \sum_\omega d(\omega)\cdot
s_\omega$.

Now, our first remark about totally resonant series is;

\proclaim (5.8.3) Proposition. (a) Every totally resonant $\Gamma$-series
is terminating, i.e., it contains only finitely many  nonzero terms
and is therefore a polynomial on $G$ in the sense of
(1.1). \hfill\break
(b) This polynomial is, moreover, identically zero unless $q(r(s))\leq 0$.

\noindent For the torus case, terminating $\Gamma$-series were
studied in [42].

\vskip .15cm

\noindent {\sl Proof:} Every representation $\rho_\omega: H\to GL(V_\omega)$
is homogeneous of degree 1 with respect to ${\bf C}^*\i H$.
Thus, in order that $R_{[\alpha]} = \bigotimes \Sigma^{\alpha(\omega)}
(V_\omega)$ satisfy $R_{[\alpha]}^H\neq 0$, it is necessary that
$\sum |\alpha(\omega)| = 0$. Let $s\in {\bf Z}^A$. Then
 $|\alpha(\omega)+s_\omega|=
|\alpha(\omega)| + d(\omega) s_\omega$. Notice also that the
denominator in the
coefficients of the series $\Phi_s(a)$ has a pole unless $\alpha(\omega)
+s_\omega \geq 0$ for all $\omega$. So our statement follows from the
next obvious fact.

\proclaim (5.8.4) Lemma. The set of $\alpha = (\alpha(\omega))_{\omega\in A}$
such that all $\alpha(\omega)\geq 0$ and $\sum |\alpha(\omega)| = d$,
is finite for all $d$ and empty for $d<0$.

\vskip .1cm

Let now $(H_0, A_0)$ be an inhomogeneous pair, and $(H,A)$ be its
homogeneization, so that we have a surjective homomorphism $p: H_0\times
{\bf C}^*\to H$ with finite kernel. Let $s\in {\bf Z}^A$
 and $\sigma=r(s)\in\Lambda^0$ be as before. Let 
$$(\sigma_0, -\tau)\in
{\rm Hom}(H_0, {\bf C}^*)\times {\rm Hom}({\bf C}^*, {\bf C}^*) = 
{\rm Hom}(H_0, {\bf C}^*)\times {\bf Z}$$
be the pullback of $\sigma$ via $p$. Clearly, $\tau = -q(r(s)) \geq 0$
if $\Phi_s$ does not vanish identically. 

\proclaim (5.8.5) Proposition. The terminating series $\Phi_s(a)$ is equal 
to the following Euler integral of $f(x) = \sum_{\omega\in A_0} (a_\omega, 
x^\omega)$:
$$\Phi_s(a) = {\rm const} \cdot \int_{H_{0,c}} f(x)^\tau x^{\sigma_0}
d^*x,$$
where $H_{0,c}\i H_0$ is a maximal compact subgroup. 

Of course, this integral can be calculated purely algebraically,
since $\tau\in {\bf Z}_+$.

\noindent {\sl Proof:} 
 We repeat the arguments of Remark 5.4.5
with $\zeta_0 = H_{0,c}$  but with ${\bf R}_+$ replaced by the
unit circle $U_1$. Accordingly, $p(\zeta)$ is replaced
by the compact form $H_c$:
$$\int_{H_{0,c}} f(x)^\tau x^{\sigma_0} dx\quad  =\quad =
{\tau !\over 2\pi i} \int_{H_{0,c}}
\left( \int_{|\lambda|=1} e^{\lambda f(x)} \lambda^{-\tau} {d\lambda
\over\lambda}\right)
x^{\sigma_0} d^*x =$$
$$= {\tau !} \int_{H_{0,c}\times U_1}
 e^{\lambda f(x)}\lambda^{-\tau} x^{\sigma_0}d^*x d^*\lambda \quad = \quad 
\tau !
D\cdot \int_{H_c} e^{\bar f (\bar x)} \bar x^\sigma d^*\bar x,$$
where $D$ is the degree of the map $p$, and $\bar f\in M_A$ corresponds
to $f\in M_{A_0}$ under the canonical identification. 
On the other hand,
$$\int_{H_c} e^{\bar f (\bar x)} \bar x^\sigma d^*\bar x \quad = \quad
\int_{G_c} \exp \biggl(\sum_\omega {\rm tr}(a_\omega b_\omega)\biggr)
\delta_H(b) b^s db.$$
Expanding now $\delta_H$ into a series of distributions on $G_c$,
as in (3.4.2), we see that this time it is legitimate to integrate it termwise,
since the integration is taken over $G_c$. The proof is finished by
applying the formula (4.7.2) to each factor of $G_c = \prod_{\omega\in A}
U_{d(\omega)}$.

\vskip .3cm

\noindent {\bf (5.9) Deformation of vanishing $\Gamma$-series.}
Let $s\in {\bf Z}^A$ be such that $\tau = -\sum d(\omega) s_\omega < 0$.
Then $\Phi_s(a)$ vanishes identically. On the Euler integral side,
however, the corresponding $A$-hypergeometric system is satisfied by functions
$J(f) = \int f(x)^{\tau} x^{\sigma_0}d^*x$, which are periods
of rational differential forms on the spherical variety $X_A$, or
on the hypersurface $\{f=0\}\i X_A$. A way to obtain solutions
from the general power series construction is by taking iterated
derivatives of the $\Phi_s(a)$ in $s$. For the toric case this method
was used by J. Stienstra [46] (see also [26]) 
to construct all the solutions
in the resonant situation. Here we will consider the simplest deformation
$s\mapsto (s_\omega + t)$, $t\in {\bf C}$. 

The function $1/\Gamma(z)$ has a first order zero at $z=-d$, $d=1,2,...$
with the derivative $(-1)^d d!$. Thus,
using the definition of the matrix $\Gamma$-function as a product, we find
$${d\over dt}\biggr|_{t=0} \Phi_{(s_\omega+t)} (a)  = \sum_{\gamma\in A}
\Phi_s^{(\gamma)}(a), \quad {\rm where} \leqno (5.9.1)$$
$$\Phi_s^{(\gamma)} (a) = \sum_{\matrix{\alpha\in\Irr :\cr
\alpha(\omega)+s_\omega \geq 0, \,\, \omega\neq \gamma, \cr
\alpha(\gamma)_{d(\gamma)-1}+s_\gamma+2 \geq 0\cr}}
\bigl( d(\alpha) C_s^{(\gamma)}(\alpha) I_\alpha(H), \,\,
a^{\alpha+s}\bigr),\leqno (5.9.1')$$
where $C_s^{(\gamma)}(\alpha)$ is the following number:
$$C_s^{(\gamma)}(\alpha) ={
(-1)^{\alpha(\gamma)_{d(\gamma)}-s_\gamma - 1}
 \Gamma \bigl(-\alpha(\gamma)_{d(\gamma)}-s_\gamma\bigr)\over
\biggl( \prod_{\omega\neq\gamma} \Gamma_{d(\omega)}(\alpha(\omega)
+s_\omega + 1)\biggr) \prod_{j=1}^{d(\gamma)-1}
\Gamma\bigl( \alpha(\gamma)_j+s_\gamma + d(\gamma)-j+1\bigr) }.
\leqno (5.9.1'')$$
(I am grateful to the referee for correcting this formula.) 
As before, we prove:

\proclaim (5.9.2) Proposition. Each $\Phi_s^{(\gamma)}$ formally
satisfies the $A$-hypergeometric system. 

We now consider one particular case, generalizing an observation of Batyrev
([1], Proposition 14.6).
Namely, let $(H_0, A_0)$ be an inhomogeneous pair such that $0\in A_0$
and, moreover, $0$ is a vertex of $Q_{A_0}$. Thus any $f\in M_{A_0}$
is written as
$$f(x) = a_0 + \sum_{\omega\in A'_0} (a_\omega, x^\omega),
\quad A'_0 = A_0-\{0\}.$$
Let $\Irr'_+ = \prod_{\omega\in A'_0} \Irr_+(GL(V_\omega))$.
For $\alpha\in\Irr'_+$ let $\hat\alpha\in\Irr$ have the following components:
$\hat\alpha(0) = -\sum_{\omega\in A'_0} d(\omega)\cdot |\alpha (\omega)|$,
while $\hat\alpha(\omega) = \alpha(\omega)$ for $\omega\in A'_0$. 
Let $(H,A)$ be the homogeneization of $(H_0, A_0)$ with the standard
bijection $\omega\mapsto \bar\omega$,  $A_0\to A$. 
Let $A'\i A$ be the image of $A'_0$. 
Taking $s_{\bar 0}= 1$ and $s_\omega=0$ for $\omega\in A'$, we get
 $s\in {\bf Z}^A$ such that $\Phi_s(a)$ vanishes identically, while
$$\Phi_s^{(0)}(a) = \sum_{\alpha\in\Irr'_+} \left(
d(\alpha) I_{\hat \alpha}(H) { \Gamma \left( 1+\sum_{\omega\in A'_0}
 |\alpha(\omega)| d(\omega) \right)\over
\prod_{\omega\in A'_0} \Gamma_{d(\omega)}(\alpha(\omega)+1)}, \,\,\,
a^{\hat\alpha}\right). \leqno (5.9.3)$$

\proclaim (5.9.4) Proposition. There is a domain $U\i M_A$ such that
for $f\in U$ the hypersurface $\{f=0\}$ does not meet the subgroup
$H_{0,c}\i H_0$, the series $\Phi_s^{(0)}(a)$ converges in $U$,
and its sum is equal to the  following 
Euler integral:
$$\Phi_s^{(0)}(a) = {\rm const} \cdot \int_{H_{0,c}} {d^*x\over f(x)}, \quad f
=\sum_{\omega\in A_0} (a_\omega, x^\omega).$$

\noindent {\sl Proof:} This is achieved, similarly to [1], by expanding
$${a_0\over f(x)} = {1 \over 1+\sum_{\omega\in A'_0} \biggl( (a_\omega/a_0),
\,\, x^\omega\biggr)}$$
into the geometric series. Terms of this geometric series are labelled by
integer vectors $m = (m_\omega)\in ({\bf Z}_+)^{A'_0}$, the term
corresponding to $m$ being
$${(-1)^{\Sigma m _\omega} \biggl(\sum m_\omega\biggr)!\over
\prod(m_\omega !)} \prod_{\omega\in A'_0} \biggl( {\rm tr} \bigl( a_\omega
x^\omega/a_0\bigr) \biggr)^{m_\omega}.\leqno (5.9.5)$$
When $a_0$  dominates all the other $a_\omega$, the geometric
series converges on a domain  containing $H_{0,c}$ and we can integrate it
termwise. 
By applying the formulas (4.6.2-4) to each factor in (5.9.5)
and using the orthogonality relations (3.0.1),
we get the identification of the integrated series with $\Phi_s^{(0)}(a)$.

\vfill\eject 

\centerline {\bf \S 6. Examples.}

\vskip 1cm

\noindent {\bf (6.1) A generalization of the Gauss function.}
We take the group $H_0 = {\bf C}^*\times GL_n = GL(L)\times GL(V)$,
$\dim(L)=1$, $\dim(V)=n$. Take $A_0 = \{ {\bf C}, L, V, L\otimes V\}$. 
The pair $(H_0, A_0)$ is inhomogeneous (2.1.4). A natural homogeneization 
of $(H_0, A_0)$ is:
$$ H =  {\bf C}^*\times {\bf C}^*\times GL_n = GL(N)\times GL(L)\times GL(V),
\quad  A = \bigl\{ N^{\otimes 2}, N\otimes L, N\otimes V, L\otimes V\bigr\},
\leqno (6.1.1)$$
where $\dim(N)=1$. Thus
$$M_{\bar A} = {\bf C}\oplus {\bf C}\oplus \Mat\oplus\Mat = 
 \End (N^{\otimes 2})\oplus \End(N\otimes L)\oplus \End(N\otimes V)
\oplus \End(L\otimes V).\leqno (6.1.2)$$
A typical point in $M_A$ will be denoted by $(a,b, C,D)$ where
$a,b\in {\bf C}$ and $C,D\in\Mat$. A point in the dual space $M_A^*$
will be denoted by $(u,v,x,y)$ where $u,v\in {\bf C}$ and $x,y\in\Mat$. 
The variety $Y_A\i M_A^*$ from (2.2) is:
$$Y_A = \{ (u,v,x,y): uy=vx\},\leqno (6.1.3)$$
and $X_A$ is the projectivization of this variety.
Thus 
$$\dim(X_A)=n^2+1, \quad {\rm deg}(X_A) = 2^{n^2}.$$
The lattice of weights of $H$ is $\Lambda = {\bf Z}\oplus {\bf Z}\oplus
{\bf Z}^n$, with the Weyl group $W=S_n$ acting by permutations
on the last summand only. We denote the basis vectors of $\Lambda$ by
$e_{-1}, e_0, e_1, ..., e_n$, with $e_{-1}$ and $e_0$ spanning the
first two summands. The polytope $Q_A$ is the cylinder over a simplex:
$$Q_A = {\rm Conv} \biggl\{ 2e_{-1}, e_{-1}+e_0, e_{-1}+e_i, e_0+e_i,
 \,\, i=1, ..., n\biggr\}\simeq \Delta^1\times\Delta^n.\leqno (6.1.4)$$
The $A$-hypergeometric system is  a system on a function $\Phi = \Phi(f) 
=\Phi(a,b,C,D)$ on $M_A$. Writing the quasihomogeneity conditions in the 
integrated from, we
can represent the system as follows:
$$\cases{ {\partial\over\partial a} {\partial\over\partial D_{ij}}\Phi = 
{\partial
\over\partial b}{\partial\over\partial C_{ij}}\Phi\cr
\Phi(\lambda^2a, \lambda b, \lambda C, D) = \lambda^{\chi_1} \Phi(a,b,C,D)\cr
\Phi(a,\lambda b, C, \lambda D) = \lambda^{\chi_2} \Phi(a,b,C,D)\cr
\Phi(a,b, UC, UD) = \Phi(a,b,CU, DU) = |U|^{\chi_3} \Phi(a,b,C,D).\cr}\leqno 
(6.1.5)$$
Here $\chi_i, i=1,2,3,$ are complex parameters labelling a character $\chi$ of
$\bar{\bf h} = {\bf C}\oplus {\bf C}\oplus gl_n$. If $\Phi$ is a solution, 
then by setting
$\phi(D)=\Phi(1,1,{\bf 1}, D)$ we have:
$$\Phi(a,b,C,D) = a^{{\chi_1\over 2}-{\chi_2\over 2} -{n\chi_3\over 2}}
 b^{\chi_2} |C|^{\chi_3}
\phi(ab^{-1}C^{-1}D).\leqno (6.1.6)$$
$$\phi(UDU^{-1})=\phi(U),\leqno (6.1.7)$$
so, in particular, $\phi$ is conjugacy invariant.
The Euler integral solutions are:
$$\Phi(a,b,C,D) = \int_{(u,y)\in \zeta\i {\bf C}^*\times GL_n}
\bigl( a+bu + (C,y)+u(D,y)\bigr)^p t^q \det(y)^r du dy,
\leqno (6.1.7)$$
where
$$p={1\over 2} \left(\chi_1+\chi_2+{\chi_3\over n}\right), \quad q=\chi_2-1, 
\quad r=\chi_3-n.
\leqno (6.1.8)$$

\proclaim (6.1.9) Proposition.  The matrix $\Gamma$-series
corresponding to $s=(s_1, ..., s_4)$
has the form
$$\Phi_s(a,b,C,D) = \sum_{\nu\in\Irr(GL_n)} {
a^{|\nu|+s_1} b^{-|\nu|+s_2} |C|^{s_3} |D|^{s_4}
 d(\nu) s_\nu(CD^{-1})\over
\Gamma(|\nu|+s_1+1) \Gamma(-|\nu|+s_2+1) \Gamma_n(\nu^- +s_3+1)
\Gamma_n(\nu+s_4+1)}.$$

\noindent {\sl Proof:} The set
$\Irr: =\Irr(G)$
consists of $\sigma = (m_1, m_2, \mu, \nu)$ where $m_i\in {\bf Z}$ and $\mu, 
\nu\in\Irr(GL_n)$.
The representation of $G$ corresponding to such a $\sigma$, has the form
$$R_\sigma = (N^{\otimes 2})^{\otimes m_1} \otimes (N\otimes L)^{\otimes m_2}
\otimes \Sigma^\mu(N\otimes V) \otimes \Sigma^\nu(L\otimes V) =$$
$$= N^{\otimes (2m_1+m_2+|\mu|)} \otimes L^{\otimes (m_2+|\nu|)} \otimes 
\Sigma^\mu (V)
\otimes\Sigma^\nu(V).$$

\proclaim (6.1.10) Proposition. (a) The space $R_\sigma$ possesses nonzero
$ H$-invariants if and only if
$\mu=\nu^-$, $m_2 = -|\nu|, m_1 = |\nu|$. In this case the space of invariants 
is 1-dimensional.
\hfill\break
(b) If the conditions of (a) are satisfied, then for $x=(a,b,C,D)\in G$ the 
function
$(I_\sigma(H), x)$ is equal to
$$(I_\sigma (H), x) ={1\over d(\nu)} a^{|\nu|} b^{-|\nu|} s_\nu (CD^{-1}).$$

\noindent {\sl Proof:} (a) is clear. 
To see (b), consider the following less cumbersome case to which everything is 
easily reduced.
Let $H' = GL(V), G' = GL(V)\times GL(V)$, so that $H'$ is the diagonal in $G'$.
Let $\sigma =(\nu, \nu^-)\in\Irr(G')$ and let $l\in R_\sigma = 
\Sigma^\nu(V)\otimes \Sigma^{\nu^-}(V) = \Sigma^\nu(V) \otimes 
(\Sigma^\nu(V))^*$
be the canonical pairing which generates the 1-dimensional space of 
$H'$-invariants.
Let $p$ be  the projection onto this 1-dimensional space along other isotypical
components, so that for $x=(C,D)\in G_0$ we have $(I_\sigma(H'), x) = {\rm tr} 
(p\circ x^\sigma)$.
This last trace can be just seen as the matrix element $\langle l | x^\sigma 
|l \rangle$
with respect to any $H'$-invariant scalar product $\langle -|-\rangle$ on 
$\Sigma^\nu(V)
\otimes \Sigma^{\nu}(V^*)$ such that $\langle l|l\rangle = 1$.  

\proclaim (6.1.11) Lemma. In the above assumptions we have, for
 $x=(C,D)\in G'$:
$$\langle l| x^\sigma |l\rangle :=
\langle l | \Sigma^\nu (C)\otimes \Sigma^{\nu^-}(D)|l\rangle =
{1\over d(\alpha)} {\rm tr}\, \Sigma^\nu(CD^{-1})
={1\over d(\alpha)} s_\nu(CD^{-1}).$$

\noindent {\sl Proof:} It is enough to consider the case $\nu = (1)$: in the 
general
case we just consider the space $\Sigma^\nu(V)$ as the new $V$ and 
$\Sigma^\nu(C)$ as the new $C$ etc.
Assuming that $\alpha = (1)$, we find that the element $l$ can be viewed as the
identity matrix ${\bf 1}\in \End(V) = V\otimes V^*$. The  $GL(V)$-invariant 
scalar product on
${\rm End}(V)$ such that $\langle {\bf 1}|{\bf 1}\rangle = 1$ is given by
$\langle X|Y \rangle = {\rm tr}(XY)/ \dim(V)$. The operator 
$\rho_\nu(C)\otimes {\bf 1}$
on $V\otimes V^*$ corresponds, after the identification $V\otimes V^*\simeq 
{\rm End}(V)$,
to the left multiplication by $C$; the operator ${\bf 1}\otimes\rho_{\nu^-}(D)$
corresponds to the right multiplication by $D^{-1}$. Thus our matrix element
is equal to
$${1\over \dim(V)} {\rm tr} ( (C{\bf 1}D^{-1})\cdot {\bf 1}) = {1\over\dim(V)} 
{\rm tr}(CD^{-1}).$$
This establishes the particular case $\nu=(1)$ and the general case of the 
lemma follows from that.
This implies, in its turn, Propositions 6.1.10 and 6.1.9. 

\vskip .1cm

The set $A$ possesses two one-element cobases (5.7.1), namely
$B_1=\{L\otimes V\}$ and $B_2=\{N\otimes V\}$. The matrix $\Gamma$
series adapted to $B_1$ can be obtained by putting $s_3=0$. 
The series $\phi_{s_1, s_2, s_4}(D) = \Phi_{s_1, s_2, 0, s_4}(1,1, {\bf 1}, D)$
has the form
$$\phi_{s_1, s_2, s_4}(D) = $$
$$=|D|^{s_4}\sum_{\mu\in \Irr^+(GL_n)} { d(\mu) s_\mu(D)\over 
\Gamma(-|\mu|+s_1+1) \Gamma(|\mu|+s_2+1) \Gamma_n(\mu+1)
\Gamma_n(\mu^-+s_4+1)}.
\leqno (6.1.12)$$
It converges near $D=0$. Introducing the standard
 Pochhammer symbol
$$(a)_m = {\Gamma (a+m)\over \Gamma(a)} = a(a+1)...(a+m-1), \quad a\in {\bf 
C}, m\in {\bf Z}_+
\leqno (6.1.13)$$
and the generalized, or matrix,  Pochhammer symbol (cf. [23])
$$[a]_\mu = \prod_{j=1}^n (a+n-j)_{\mu_j} = {\Gamma _n(\mu_1 + a, ..., 
\mu_n+a)\over
\Gamma_n(a, ..., a)}, 
\quad a\in {\bf C}, \mu\in \Irr^+(GL_n) \leqno (6.1.14)$$
we can write, by using (4.6.4):
$$\phi_{s_1, s_2, s_4}(D) = {\rm const}\cdot |D|^{s_4} \cdot \, _2{\cal F}_1 
(-s_1, -s_4,
s_2+1; D),$$
where 
$$_2{\cal F}_1(\alpha, \beta, \gamma; x) = \sum_{\mu\in\Irr^+(GL_n)}
{ (\alpha)_{|\mu|} [\beta]_\mu\over (\gamma)_{|\mu|} \, |\mu|!}
w_\mu s_\mu(x), \quad x\in\Mat. \leqno (6.1.15)$$
This function is similar to but not identical with the James-Biedenharn-Louck
matrix
generalization of the Gauss function [24] [30] [36], which is
$$_2F_1(\alpha,\beta, \gamma; x) = \sum_{\mu\in\Irr^+(GL_n)} 
{[\alpha]_\mu \cdot [\beta]_\mu \over [\gamma]_\mu \cdot |\mu|!}  w_\mu 
s_\mu(x),
\quad x\in\Mat, \|x\|< 1. \leqno (6.1.16)$$

\vskip .3 cm

\noindent {\bf (6.2) A generalization of the Pochhammer function.} 
Set $G= ({\bf C}^*)^{2p}\times GL_n\times GL_n$ and let
$$H=\biggl\{ (u_1, ..., u_p, v_1, ..., v_p, x,y)\in G\biggl|
u_1...u_p\cdot y = v_1...v_p\cdot x\biggr\}.\leqno (6.2.1)$$
Here $x,y\in GL_n$. Thus $H\simeq ({\bf C}^*)^{2p}\times GL_n$. Let $A$ be the
set of the tautological  representations of the factors of $G$
considered as representations of $H$. Thus $Y_A\i M_A^* =
{\bf C}^{2p}\times\Mat\times\Mat$ is given by the same equation as
in (6.2.1). A point of $M_A$ will be denoted by
$(a_1, ..., a_p, b_1, ..., b_p, C,D)$. The subgroup $H_0\i H$ given by
$u_1=1$ satisfies the inhomogeneity condition. The $A$-hypergeometric
system depends on $2p+1$ complex parameters $\chi_1, ..., \chi_{2p+1}$
and consists of the equations
$${\partial^{p+1}\over\partial a_1 ... \partial a_p \partial d_{ij}}\Phi = 
{\partial^{p+1}\over\partial b_1 ... \partial b_p \partial c_{ij}}\Phi$$
together with the quasihomogeneity conditions involving the $\chi_i$ which
we leave to the reader. The Euler integral solutions have the form
$$\Phi(a_1, ..., a_p, b_1, ..., b_p, C,D) =$$
$$= \int_\zeta \biggl(a_0 + \sum_{i=2}^p a_iu_i + \sum_{i=1}^p b_iv_i + (C,y) +
{v_1...v_n\over u_2, ..., u_n}(D,y)\biggr)^{\lambda_1}\times$$
$$\times \prod_{i=2}^p u_i^{\lambda_i}\prod_{i=1}^p v_i^{\lambda_{p+i}}
\det(y)^{\lambda_{2p+1}}dy \prod_{i=2}^p du_i \prod_{i=1}^p dv_i,$$
where $\zeta$ is an appropriate $2p-1+n^2$-dimensional contour in $H_0$. 
Any of the two $n$-dimensional representations in $A$ forms a cobase
with one element. The corresponding convergent $\Gamma$-series
are easily expressed through the function
$$_{p+1}{\cal F}_p (\alpha_1, ..., \alpha_{p+1}, \beta_1, ..., \beta_p; x) = 
\sum_{\mu\in\Irr^+(GL_n)} { (\alpha_1)_{|\mu|} ... (\alpha_p)_{|\mu|}
[\alpha_{p+1}]_\mu\over (\beta_1)_{|\mu|} ... (\beta_p)_{|\mu|} \,
|\mu|!} w_\mu s_\mu(x).$$

\vskip .2cm

\noindent {\bf (6.3) A generalization of the Appell's  functions
via $3j$-symbols.} We consider the same groups $H_0$ and $H$ as in (6.1)
but take a bigger set of representations of $H_0$, namely
$$A_0 = \bigl\{ {\bf C}, L, V, L\otimes V, L^{\otimes 2}\otimes V\bigr\} \i 
\Irr(H_0)\leqno (6.3.1)$$
which entails for the homogeneization:
$$A=\biggl\{ N^{\otimes 3}, N^{\otimes 2}\otimes L, N^{\otimes 2}\otimes V,
 N\otimes L\otimes V, L^{\otimes 2}\otimes V\biggr\}.\leqno (6.3.2)$$
The space $M_A$ consists of $(a,b, x,y,z)$ with $a,b\in {\bf C}$ and 
$x,y,z\in \Mat$. 
 
We will concentrate here on the form of
matrix $\Gamma$-series. 
The group $G$ is ${\bf C}^*\times {\bf C}^* \times GL_n\times GL_n\times GL_n$.
Irreducible  representations of $G$ are parametrized by tuples
 $\sigma =(m,r, \lambda,\mu, \nu)$ with $m,r\in {\bf Z}$ and $\lambda, 
\mu, \nu\in\Irr(GL_n)$. The corresponding representation is
$$R_\sigma = N^{\otimes (3m+2n + 2|\lambda| + |\mu|)} \otimes L^{\otimes
(r + |\mu| + 2|\nu|)} \otimes \Sigma^\lambda(V)\otimes \Sigma^\mu V \otimes
\Sigma^\nu V.$$
The problem of finding $H$-invariants in $R_\sigma$ is thus identical
with the Clebsch-Gordan problem for the group $GL_n$. When $n=1$,
the problem is trivial, and among the hypergeometric series one
finds Appell's functions $G_1, H_3, H_6$, as explained in [15], \S 3.2. 
So we consider the case $n=2$. In this case the elementary
representation theory of $SL_2$ implies  the following.

\proclaim (6.3.3) Proposition. (a) The necessary and sufficient conditions
for non-vanishing of $R_\sigma^H$ are:\hfill\break
(a1) $3m+2r+2|\lambda| + |\mu| = r+|\mu|+2|\nu| = 0$,
\hfill\break
(a2)  $|\lambda+ |\mu| + |\nu| = 0$ and the nonnegative integers 
$\lambda_1-\lambda_2, \,\,\mu_1-\mu_2, \,\, \nu_1-\nu_2$ satisfy the
triangle inequalities, i.e., each of them does not exceed the sum of
the two others.\hfill\break
(b) If the conditions of (a) are satisfied, 
then $R_\sigma^H$ is one-dimensional.

To actually write down the  series, we need to fix a basis
in $\Sigma^\lambda V$, $V={\bf C}^2$. We choose the Gelfand-Cetlin
convention, so that the basis vectors are denoted by $e^{(\lambda)}_k$
with $k$ an integer, $\lambda_1\geq k\geq \lambda_2$. In this basis,
the matrix elements are the following functions of $x=\|x_{ij}\|\in GL_2$
(see [48]. p. 116):
$$t^{(\lambda)}_{km}(x) = |x|^{\lambda_2}
\sqrt{ (k-\lambda_2)! (\lambda_1-k)!\over (m-\lambda_2)! (\lambda_1-m)!}
\sum_{i+j=m-\lambda_2\atop i,j\geq 0} {k-\lambda_2\choose i}
{\lambda_1-k \choose j} x_{11}^i x_{12}^{k-\lambda_2-i} x_{21}^j x_{22}^
{\lambda_1-k-j}.$$
We will use the hermitian scalar product in which the $e^{(\lambda)}_k$ form
 an orthonormal basis. It has the property that the action of $U_2\i GL_2$
is unitary. If $\lambda, \mu, \nu$ satisfy the conditions of (6.3.3)(a2),
then the normalized 
generator of the 1-dimensional $GL(V)$-invariant subspace in
$\Sigma^\lambda V \otimes\Sigma^\mu V\otimes\Sigma^\nu V$ is
$$v_{\lambda\mu\nu} = \sum_{i,j,k} \pmatrix{\lambda&\mu&\nu\cr i&j&k}
e^{(\lambda)}_i\otimes e^{(\mu)}_j\otimes e^{(\nu)}_k, \quad
\| v_{\lambda  \mu\nu}\|=1,$$
where the coefficients, called the 3j-symbols, are
nonzero only for $i+j+k=0$. They are easily reduced to
the Clebsch-Gordan coefficients for the group $SU_2$, see [19] [37] [48]  and 
are, explicitly, as follows ([37], \S 54):
$$\pmatrix{\lambda&\mu&\nu\cr i&j&k} = 
\sqrt{ (\lambda_2+\mu_1+\nu_1)! (\lambda_1+\mu_2+\nu_1)!\over
 (2\lambda_1+2\mu_1+2\nu_1+1)!}\times$$
$$\times \sqrt{ (j-\mu_2)! (\mu_1-j)! (k-\nu_2)!
(\nu_1-k)! (\lambda_1+\mu_1+
\nu_2)! (i-\lambda_1)! (\lambda_1-i)!}\times$$
$$\times \sum_z {(-1)^{z+\nu_1-k}\over
z! (\lambda_1+\mu_1+\nu_2-z)! (z-j-\nu_2-\lambda_1)! (j-z-\mu_2)!
(z+i+\lambda_1+\mu_2+\nu_1)! (\lambda_1-z-i)!}.$$
The matrix $\Gamma$-series corresponding to $s=(s_1, ..., s_5)$ has the 
form:
$$\Phi_s(a,b,x,y,z) = a^{s_1} b^{s_2} |x|^{s_3} |y|^{s_4} |z|^{s_5} \cdot
\sum_{\lambda,\mu, \nu\in\Irr(GL_2)\atop |\lambda|+|\mu|+|\nu|=0}\left(
d(\lambda)d(\mu)d(\nu)({b\over a})^{|\mu|+2|\nu|}\times\right. $$
$$\left.\times {  \sum_{i,j,k,i',j',k'}
\pmatrix{\lambda&\mu&\nu\cr i&j&k}\pmatrix{\lambda&\mu&\nu\cr i'&j'&k'}
t^{(\lambda)}_{ii'}(x) t^{(\mu)}_{jj'}(y) t^{(\nu)}_{kk'}(z)\over
\Gamma(-|\mu|-2|\nu|+s_1+1)\Gamma (|\mu|+2|\nu|+s_2+1)
 \Gamma_2(\mu+s_3+1)\Gamma_2(\nu+s_4+1)
\Gamma_2(\nu+s_5+1)}\right).$$
Any two 2-dimensional representations from $A$ form a cobase. The three
types of convergent $\Gamma$ series  adapted to these
cobases, are natural $GL_2$-generalizations of Appell's series
$G_1, H_3, H_6$.

\vskip 2cm   

\centerline {\bf References.}

\vskip 1cm

\item{[1]} V.V. Batyrev, Variations of the mixed Hodge structure of affine
hypersurfaces in algebraic tori, {\it Duke Math. J.} {\bf 69} (1993),
349-409. 

\item{[2]} A.Borel et al, Algebraic D-modules, Academic Press, 1987. 

\item{[3]} M. Brion, Spherical varieties: an introduction, in: ``Topological 
methods
in Algebraic Transformation Groups" (H. Kraft et al. Eds.) (Progr. Math {\bf 
80}),
p. 11-26, Birkhauser, Boston, 1989.

\item{[4]} M. Brion, Groupe de Picard et nombres characteristiques des 
vari\'etes
sph\'eriques, {\it Duke Math. J.} {\bf 58} (1989), 397-424.

\item{[5]} M. Brion, D. Luna, T. Vust, Espaces homog\`enes sph\'eriques,
{\it Invent. math.} {\bf 84} (1986), 617-632. 

\item{[6]}J.-L.Brylinski, Transformations canoniques, dualit\'e projective,
th\'eorie de Lefschetz, transformations de Fourier et sommes
trigonom\'etriques, {\it Ast\'erisque} {\bf 140-141} (1986), p. 3-134.

\item{[7]} C. De Concini, C. Procesi, Complete symmetric varieties,
in: ``Invariant Theory" (F. Gherardelli, Ed.)
Lecture Notes in Math. {\bf 996}, p. 1-44.

\item{[8]} C. De Concini, C. Procesi, Complete symmetric varieties II, 
in: ``Algebraic Groups and Related Topics" (Adv. Studies in Pure Math. {\bf 6})
p. 481-513, Academic Press, 1985.

\item{[9]} A. Erdelyi, W. Magnus, F. Oberhettinger, G.F. Tricomi,
 Higher Transcendental Functions, vol. 1,
Mc Graw-Hill Publ. New York, 1953.

\item{[10]} W. Fulton, Introduction in the Theory of Toric Varieties, 
Princeton Univ. Press, 1995.

\item{[11]} W. Fulton, J. Harris, Representation Theory 
(Graduate Texts in Math. 
{\bf 129}) Springer-Verlag 1991. 

\item{[12]} I.M. Gelfand, M.I. Graev, V.S. Retakh, Generalized hypergeometric
functions associated to any finite or locally compact field,
{\it Russian Math. Dokl.} {\bf 46} (1993), 343-347.

\item{[13]}  I.M. Gelfand, M.I. Graev, V.S. Retakh, Generalized hypergeometric
systems of equations
and series of hypergeometric type,
 {\it Russian Math. Surv.} {\bf 47} (1992),  N. 4, p. 1-88.

\item{[14]} I.M. Gelfand, M.I. Graev, A.V. Zelevinsky, Holonomic systems of 
equations
and series of hypergeometric type, Dokl. AN. SSSR, {\bf 295} (1987) p. 14-19.

\item{[15]} I.M. Gelfand, M.M. Kapranov, A.V. Zelevinsky, Hypergeometric 
functions
and toric varieties, Funct. Anal. Appl. {\bf 23}  (1989), p. 94-106, see
also Correction, {\it ibid.} {\bf 27} (1994), 295.

\item{[16]} I.M. Gelfand, M.M. Kapranov, A.V. Zelevinsky. Generalized Euler 
integrals and
$A$-hypergeometric functions, {\it Adv. Math.} {\bf 84} (1990), 255-271. 

\item{[17]} I.M. Gelfand,  M.M. Kapranov,  A.V. Zelevinsky. Hypergeometric 
functions,
toric varieties and Newton polyhedra, in: ``Special Functions" (M. Kashiwara, 
T. Miwa Eds.)
(ICM-90 Satellite Conference Proceedings) p. 104-121, Springer-Verlag, 1991. 

\item{[18]} I.M. Gelfand, M.M. Kapranov,  A.V. Zelevinsky, Discriminants, 
Resultants and
Multidimensional Determinants, Birkh\"auser, Boston, 1994.

\item{[19]} I.M. Gelfand, R.A. Minlos, Z.Y. Shapiro, Representations of the 
Rotation
Group and  Lorentz Group and Their Applications, Pergamon Press, New York, 
1963.

\item{[20]} I.M. Gelfand, G.E. Shilov, Generalized Functions I: Generalized
Functions and Operations on Them, Academic Press, New York 1964.

\item{[21]} K.I. Gross, R. A. Kunze, Fourier  transforms and
holomorphic discrete series, in: ``Conference on Harmonic Analysis"
(D.Gulick, R. Lipsman Eds.), Springer Lecture Notes in Math. 
{\bf 266}, p. 79-122, Springer-Verlag 1972.

\item{[22]} K.I. Gross, R.A. Kunze, Bessel functions and representation theory 
II:
Holomorphic discrete series and metaplectic representations,
{\it J. Funct. Anal.} {\bf 25} (1977), 1-49. 

\item{[23]} K.I. Gross, D.St.P. Richards, Special functions of matrix 
argument. I: Algebraic
induction, zonal polynomials and hypergeometric functions, {\it Trans. AMS},
{\bf 301} (1987), 781- 811. 

\item{[24]} K.I. Gross, D.St.P. Richards, Hypergeometric functions on complex
matrix space, {\it Bull. AMS}, {\bf 24} (1991), 349-355. 

\item{[25]} S. Hosono, A. Klemm, S. Theisen, S.-T. Yau, Mirror
symmetry, mirror map and applications to complete intersection
Calabi-Yau spaces, in:  ``Mirror Symmetry II"
(B. Greene, S.-T. Yau, Eds.), International Press, Boston, 1997.

\item{[26]} S. Hosono, B.H. Lian, S.-T. Yau, GKZ-generalized hypergeometric 
systems
in mirror symmetry of Calabi-Yau hypersurfaces, {\it Comm. Math. Phys.}
{\bf 182} (1996), 535-577

\item{[27]} S. Hosono, B.H. Lian, S.-T. Yau, Maximal degeneracy points of GKZ 
systems,
{\it J. Amer. Math. Soc.} {\bf 10} (1997), 427-443.

\item{[28]} R. Hotta, Equivariant D-modules, to appear in: `` Proc. of the
ICPAM Spring School in Wuhan",  R. Torasso Ed., Hermann. 

\item{[29]} L.K. Hua, Harmonic Analysis of Functions of Several Variables in
the Classical Domains (Translations of Math.Monographs {\bf 6}),
Amer. Math. Soc. 1963.

\item{[30]} A.T. James, Distributions of matrix variates and latent roots 
derived from
normal samples, {\it Ann. Math. Stat.} {\bf 35} (1964), 475-501.

\item{[31]} B. Ya. Kazarnovskii, Newton polyhedra and Bezout formula for
matrix-valued functions of finite-dimensional representations,
{\it Funct. Anal. Appl.} {\bf 21} (1986), 319-321. 

\item{[32]} S. Kerov, The boundary of Young lattice and random Young tableaux,
in: ``Formal Power Series and Algebraic Combinatorics" (L. Billera et al.
Eds.) (DIMACS series in discrete mathematics {\bf 24})
p. 133-158, American Math. Soc. 1996. 

\item{[33]} F. Knop, The Luna-Vust theory of spherical embeddings, in: 
Proc. Hyderabad Conf. on Algebraic Groups (S. Ramanan Ed.) p. 225-249,
Manoj Prakashan Publ., Madras, 1991.  

\item{[34]} A. Kouchnirenko, Poly\`edres de Newton et nombres de Milnor,
{\it Invent. Math.} {\bf 32} (1976), 1-31. 

\item{[35]} M. Lassalle, Sur la transformation de Fourier-Laurent dans un
groups analytique complexe r\'eductif, {\it Ann. Inst. Fourier, Grenoble},
{\bf 28} (1978), 115-138. 

\item{[36]} J.D. Louck, L.C.Biedenharn, A generalization of the Gauss 
hypergeometric
function, {\it J. Math. Anal. Appl.} {\bf 59} (1977), 423-431.

\item{[37]} G.Ja. Liubarskii, Group Theory and its Applications to Physics,
Pergamon Press, New York, 1960.

\item{[38]} H. Matsumura, Commutative Ring Theory,
Cambridge Univ. Press, 1986.

\item{[39]} D. Mumford, J. Fogarthy, F. Kirwan, Geometric Invariant Theory
(Ergebnisse der math. {\bf 34}), Springer-Verlag, 1994.

\item{[40]} F. Pham, La descente des cols par les onglets de Lefschetz,
avec vues de Gauss-Manin, {\it Ast\'erisque} {\bf 130} (1985),
11-47.

\item{[41]} A. Rittatore, Algebraic monoids and group embeddings,
preprint math.AG/9802073.

\item{[42]} M. Saito, B. Sturmfels, N. Takayama, Hypergeometric
polynomials and integer programming, {\it Compositio Math.} to appear.

\item{[43]} J.-P. Serre (avec la collaboration de P. Gabriel), Alg\`ebre 
locale. Multiplicit\'es, Lecture Notes in Math. {\bf 11}, Springer-Verlag,
1965.

\item{[44]} F.J. Servedio, Prehomogeneous vector spaces and varieties,
{\it Trans. AMS}, {\bf 176} (1973), 421-444.

\item{[45]} T. Springer, Reductive groups, in: ``Automorphic Forms,
Representations and L-functions" (A. Borel, W. Casselman Eds.)
 (Proc. Symp. Pure Math. XXXIII) 
part 1, 
p. 3-27, Amer. Math. Soc, 1979. 

\item{[46]} J. Stienstra, Resonant hypergeometric systems and mirror
symmetry,
preprint alg-geom/9711002.  

\item{[47]} B. Sturmfels, N. Takayama, Gr\"obner bases and hypergeometric 
functions,
preprint RIMS-1166, RIMS, Kyoto Univ., 1997.

\item{[48]} N.J. Vilenkin, Special Functions and the Theory of Group
Representations, AMS, Providence, 1968.  

\item{[49]} E.B. Vinberg, On reductive algebraic semigroups, in: ``Lie Groups 
and
Lie Algebras: E.B. Dynkin's seminar" (S.G. Gindikin, E.B. Vinberg, Eds.)
(Advances in the Math. Sciences, {\bf 169}), p. 145-182. 

\vskip 3cm

\noindent {\sl Author's address; Department of Mathematics,
Northwestern University, Evanston IL 60208 USA}

\vskip .2cm

\noindent {\sl email: kapranov@math.nwu.edu}

\bye